\documentclass[12pt,leqno]{article}
\usepackage{amsthm,amsfonts,amssymb,amsmath,eufrak,oldgerm}
\usepackage{epsfig}
\newcommand{\de}{\delta}
\newcommand{\des}{{\delta_*}}
\newcommand{\dde}{\dot\delta}
\renewcommand\d{\partial}

\newcommand\R{\mathbb R}

\def\eps{\varepsilon }
\newcommand\br{\begin{rem}}
\newcommand\er{\end{rem}}
\newcommand\bp{\begin{pmatrix}}
\newcommand\ep{\end{pmatrix}}
\newcommand\be{\begin{equation}}
\newcommand\ee{\end{equation}}
\newcommand\ba{\begin{equation}\begin{aligned}}
\newcommand\ea{\end{aligned}\end{equation}}

\newcommand{\CalL}{\mathcal{L}}

\newcommand{\CalR}{\mathcal{R}}
\newcommand{\CalS}{\mathcal{S}}
\newcommand{\CalT}{\mathcal{T}}

\newcommand{\ubar}{\bar u}
\newtheorem{theo}{Theorem}[section]
\newtheorem{prop}[theo]{Proposition}
\newtheorem{cor}[theo]{Corollary}
\newtheorem{lem}[theo]{Lemma}
\newtheorem{defi}[theo]{Definition}

\newtheorem{rem}[theo]{Remark}

\numberwithin{equation}{section}
\title{
Stability of undercompressive viscous shock profiles of
hyperbolic--parabolic systems.
}

\author{\sc \small Mohammadreza Raoofi\thanks{
Max Planck Institute for Mathematics in Sciences,
Leipzig, Germany;
raoofi@mis.mpg.de:
Research of M.R. was partially supported
by Max Planck Institute.
} and
Kevin Zumbrun\thanks{Indiana University, Bloomington, IN 47405;
kzumbrun@indiana.edu:
Research of K.Z. was partially supported
under NSF grant number DMS-0300487.
 }}
\begin{document}

\maketitle

\begin{abstract}
Extending to systems of hyperbolic--parabolic conservation laws
results of Howard and Zumbrun for strictly parabolic systems,
we show for viscous shock profiles of arbitrary amplitude and
type that necessary spectral (Evans function) conditions for
linearized stability established by Mascia and Zumbrun are also
sufficient for linearized and nonlinear phase-asymptotic stability,
yielding detailed pointwise estimates and sharp
rates of convergence in $L^p$, $1\le p\le \infty$.
\end{abstract}

%\clearpage
%\tableofcontents
%\clearpage
%%%%%%%%%%%%%%%%%
\bigbreak

\section{Introduction}\label{introduction}
%%%%%%%%%%%%%%%%%%%%%%%%%%%%%%%%%%%%%%%%%%%%%%%%%%%%%%%%%%%%
Consider a (possibly) large-amplitude {\it viscous shock profile},
or traveling-wave solution
\begin{equation} \label{profile}
\bar{u} (x - st); \quad \lim_{x \to \pm \infty} \bar{u} (x) =
u_{\pm},
\end{equation}
of a system of partially or fully parabolic conservation laws
\begin{equation} \label{main}
\begin{aligned}
u_t + F(u)_x &= (B(u)u_x)_x, \\
\end{aligned}
\end{equation}
$x \in \mathbb{R}$, $u$, $F \in \mathbb{R}^n$, $B\in\mathbb{R}^{n
\times n}$, where
\begin{equation}
u = \begin{pmatrix} u^I \\ u^{II} \end{pmatrix}, \quad F =
\begin{pmatrix} F^I \\ F^{II} \end{pmatrix}, \quad
 B =
\begin{pmatrix} 0 & 0 \\ b_1 & b_2 \end{pmatrix},
\end{equation}
$u^I \in \mathbb{R}^{n-r}$, $u^{II} \in \mathbb{R}^r$, $r$ some
positive integer, possibly $n$ (full regularization), and
\begin{equation*}
\text{Re }\sigma(b_2) \ge \theta > 0.
\end{equation*}
Here and elsewhere, $\sigma$ denotes spectrum of a matrix or other
linear operator. Working in a coordinate system moving along with
the shock, we may without loss of generality consider a standing
profile $\bar{u} (x)$, $s=0$; we take $s=0$ from now on.

\medskip
Following \cite{Z4}, we assume that, by some invertible change of
coordinates $u\to w(u)$, followed if necessary by multiplication
on the left by a nonsingular matrix function $S(w)$, equations
(\ref{main}) may be written in the {\it quasilinear, partially
symmetric hyperbolic-parabolic form}
\begin{equation}
\tilde A^0 w_t +  \tilde A w_{x}= (\tilde B w_{x})_{x} + G, \quad
w=\left(\begin{matrix} w^I
\\w^{II}\end{matrix}\right), \label{symm}
\end{equation}
$w^I\in \mathbb{R}^{n-r}$, $w^{II}\in \mathbb{R}^r$, $x\in
\mathbb{R}$, $t\in \mathbb{R}_+$, where, defining $w_\pm:= w(u_\pm)$:\\
\medskip
(A1)\quad $\tilde A(w_\pm)$, $\tilde A_{11}$, $\tilde
A^0$ are symmetric, $\tilde A^0 >0$.\\
\medskip
(A2)\quad Dissipativity: no eigenvector of $ dF(u_\pm)$ lies in
the kernel of $ B(u_\pm)$. (Equivalently, no eigenvector of $
\tilde A (\tilde A^0)^{-1} (w_\pm)$
lies in the kernel of $\tilde B (\tilde A^0)^{-1}(w_\pm)$.)\\
\medskip
(A3) \quad $ \tilde B= \left(\begin{matrix} 0 & 0 \\ 0 & \tilde b
\end{matrix}\right) $, $ \tilde G= \left(\begin{matrix}  0 \\ \tilde g\end{matrix}\right) $,
with $ Re  \tilde b(w) \ge \theta $ for some $\theta>0$, for all
$w$, and $\tilde g(w_x,w_x)=\mathbf{O}(|w_x|^2)$.

\medskip
Here, the coefficients of (\ref{symm}) may be expressed in terms
of the original equation (\ref{main}), the coordinate change $u\to
w(u)$, and the approximate symmetrizer $S(w)$, as
\begin{equation}
\begin{aligned} \tilde A^0&:= S(w)(\partial u/\partial w),\quad
\tilde A:= S(w)d F(u(w))(\partial u/\partial w),\\
\tilde B&:= S(w)B(u(w))(\partial u/\partial w), \quad G= -(dS
w_{x}) B(u(w))(\partial u/\partial w) w_{x}.
\end{aligned}
\label{coeffs}
\end{equation}
Alternatively, we assume, simply,

\medskip
(B1)\quad Strict parabolicity: $n=r$, or, equivalently, $\Re \sigma(B)>0$.\\
\medskip

Along with the above structural assumptions, we make the technical
hypotheses:\\
\medskip
(H0)\quad $F$, $B$, $w$, $S\in C^{5}$.\\
\smallskip
(H1)\quad The eigenvalues of $\tilde
A_*:=\tilde{A}_{11}(\tilde{A}^0_{11})^{-1}$ are (i) distinct from
$0$; (ii) of common sign; and (iii) of constant multiplicity with
respect to
$u$.\\
\smallskip
(H2)\quad The eigenvalues of $dF(u_\pm)$ are real, distinct, and
nonzero.\\
\smallskip
({H3})\quad Nearby $\bar u$, the set of all solutions of
\eqref{profile}--\eqref{main} connecting the same values $u_\pm$
forms a smooth manifold $\{\bar u^\delta\}$, $\delta\in
\mathcal{U}\subset \mathbb{R}^\ell$, $\bar u^0=\bar u$.

\begin{rem}\label{profileegs}
\textup{
Structural assumptions (A1)--(A3) [alt. (B1)]
and technical hypotheses (H0)--(H2) admit such physical systems
as the compressible Navier--Stokes equations, the equations
of magnetohydrodymics, and Slemrod's model for van der Waal gas
dynamics \cite{Z4, Z5}.   Moreover, existence of waves $\bar{u}$
satisfying (H3) has been established in each of these cases.}
\end{rem}

\bigskip
\begin{defi}\label{type}
An ideal shock
\begin{equation} \label{shock}
u(x,t) =
\begin{cases}
u_- &x < st, \\
u_+ &x > st,
\end{cases}
\end{equation}
is classified as {\it undercompressive}, {\it Lax}, or {\it
overcompressive} type according as $i-n$ is less than, equal to,
or greater than $1$, where $i$, denoting the sum of the dimensions
$i_-$ and $i_+$ of the center--unstable subspace of $dF(u_-)$ and
the center--stable subspace of $dF(u_+)$, represents the total
number of characteristics incoming to the shock.

A viscous profile \eqref{profile} is classified as {\it pure
undercompressive} type if the associated ideal shock is
undercompressive and $\ell=1$, {\it pure Lax} type if the
corresponding ideal shock is Lax type and $\ell=i-n$, and {\it
pure overcompressive} type if the corresponding ideal shock is
overcompressive and $\ell=i-n$, $\ell$ as in (H3). Otherwise it is
classified as {\it mixed under--overcompressive} type; see
\cite{ZH}.
\end{defi}

Pure Lax type profiles are the most common type,
and the only type
arising in standard gas dynamics, while pure over- and undercompressive type
profiles arise in magnetohydrodynamics (MHD) and phase-transitional models.
Mixed under--overcompressive profiles are also possible,
as described in \cite{LZ2, ZH}, 
%but we do not know a physical example.
though we do not know a physical example.
In the pure Lax or undercompressive case,
$\{\bar u^\delta\}=\{\bar u(\cdot-\delta)\}$ is just the
set of all translates of the base profile $\bar u$,
whereas in other cases it involves also deformations of $\bar u$.
For further discussion of existence, structure, and
classification of viscous profiles, see, e.g.,
\cite{LZ2, ZH, MZ2, MZ3, MZ4, Z3, Z4, Z5}, and references therein.

\begin{defi}\label{orbital}
The profile $\bar u$ is said to be {\it nonlinearly orbitally
stable} if $\tilde u(\cdot,t)$ approaches $\bar u^{\delta(t)}$ as
$t\to \infty$, $\bar u^\delta$ as defined in (H3), for any
solution $\tilde u$ of (\ref{main}) with initial data sufficiently
close in some norm to the original profile $\bar u$. If, also, the
phase $\delta(t)$ converges to a limiting value $\des$,
the profile is said to be {\it nonlinearly phase-asymptotically
orbitally stable}.
\end{defi}
\medbreak

An important result of \cite{MZ3}
was the identification of the following {\it stability criterion}
equivalent to $L^1 \to L^p$
linearized orbital stability of the profile, $p>1$,
where $D(\lambda)$ as described in \cite{GZ, ZH}
denotes the Evans function associated with the linearized
operator $L$ about the profile:
an analytic function analogous to
the characteristic polynomial of a finite-dimensional operator,
whose zeroes away from the essential spectrum agree in location
and multiplicity with the eigenvalues of $L$.

\medbreak
($\mathcal{D}$) \quad There exist
precisely $\ell$ zeroes of $D(\cdot)$ in the nonstable half-plane
$\R \lambda \ge 0$, necessarily at the origin $\lambda=0$.
\medbreak

As discussed, e.g., in \cite{MZ4, Z4, Z5}, under assumptions
(A1)--(A3) and (H0)--(H3), ($\mathcal{D}$) is equivalent to (i)
{\it strong spectral stability}, $\sigma(L)\subset \{\R \lambda <
0\}\cup \{0\}$, (ii) {\it hyperbolic stability} of the associated
ideal shock, and (iii) {\it transversality} of $\bar u$ as a
solution of the connection problem in the associated
traveling-wave ODE, where hyperbolic stability is defined for Lax
and undercompressive shocks by the Lopatinski condition of
\cite{M1, M2, M3, Fre} and for overcompressive shocks by an
analogous long-wave stability condition \cite{Z3, Z4}. Here and
elsewhere $\sigma$ denotes spectrum of a linearized operator or
matrix.

The stability condition holds always for small-amplitude Lax profiles
%\cite{Go.1}, \cite{MN}, \cite{KM}, \cite{KMN},
\cite{HuZ, PZ, FreS},
but may fail for large-amplitude, or nonclassical over- or undercompressive
profiles \cite{AMPZ1, GZ, FreZ, ZS, Z4}.
It may be readily checked numerically, as described, e.g., in
\cite{Br1, Br2, Br3, BrZ, BDG}.
It was shown by various techniques
in \cite{MZ1, MZ2, MZ3, MZ4, MZ5, Ra, HR, HRZ}
that the linearized stability condition ($\mathcal{D}$) is also
sufficient for nonlinear orbital stability of Lax or overcompressive profiles
of arbitrary amplitude.
In the strictly parabolic case (B1), this result was
extended in \cite{HZ} to shocks of arbitrary amplitude {\it and type},
in particular to shocks of under- or mixed over--undercompressive type.
However, up to now, it had not been verified for under- or
under--overcompressive profiles of systems with real viscosity.

In this paper, we establish for shocks of any type and general systems
satisfying (A1)--(A3) [alt. (B1)] and (H0)--(H3) that
($\mathcal{D}$) is sufficient for nonlinear phase-asymptotic
orbital stability.
More precisely, denoting by
\begin{equation}\label{aj}
a_1^\pm<a_2^\pm < \cdots < a_n^\pm
\end{equation}
the eigenvalues of the limiting convection matrices $A_\pm:= df(u_\pm)$,
define
\begin{equation}\label{theta}
\theta(x,t):=
\sum_{a_j^-<0}(1+t)^{-1/2}e^{-|x-a_j^-t|^2/Lt}
+ \sum_{a_j^+>0}(1+t)^{-1/2}e^{-|x-a_j^+t|^2/Lt},
\end{equation}
\begin{equation}\label{psi1}
\begin{aligned}
\psi_1(x,t)&:=
\chi(x,t)\sum_{a_j^-<0}
(1+|x|+t)^{-1/2} (1+|x-a_j^-t|)^{-1/2}\\
&\quad+
\chi(x,t)\sum_{a_j^+>0}
(1+|x|+t)^{-1/2} (1+|x-a_j^+t|)^{-1/2},\\
\end{aligned}
\end{equation}
and
\begin{equation}\label{psi2}
\begin{aligned}
\psi_2(x,t)&:=
(1-\chi(x,t)) (1+|x-a_1^-t|+t^{1/2})^{-3/2}\\
&\quad +(1-\chi(x,t)) (1+|x-a_n^+t|+t^{1/2})^{-3/2},
\end{aligned}
\end{equation}
where $\chi(x,t)=1$ for $x\in [a_1^-t, a_n^+t]$ and zero otherwise,
and $L>0$ is a sufficiently large constant.

Then, we have the following main theorem.
%NOTE: generalizing Theorem TODO of \cite{HZ} in the strictly parabolic case.

\begin{theo}\label{nonlin}
 Under assumptions
(A1)--(A3) [alt. (B1)], (H0)--(H3), and $(\mathcal D)$, the
profile $\bar u$ is nonlinearly phase-asymptotically orbitally
stable with respect to $H^{4}$
 initial perturbations $u_0$,
with $\|(1+|x|^2)^{-3/4}u_0(x)\|_{H^5}\le E_0$ sufficiently small.
More precisely, there
%ENDCHANGED
exist $\delta(\cdot)$ and $\des$ such that
\begin{equation}
\label{pointwise}
\begin{aligned}
|\tilde u(x,t)-\bar u^{\des + \delta(t)}(x)|&\le C E_0
(\theta+\psi_1+\psi_2)(x,t),\\
|\partial_x\big(\tilde u(x,t)-\bar u^{\des + \delta(t)}(x)\big)|&\le C E_0
(\theta+\psi_1+\psi_2)(x,t),\\
|\des|&\le CE_0,\\
  |\dot \delta (t)|&\le C E_0 (1+t)^{-1},\\
  |\delta(t)|&\le C E_0(1+t)^{-1/2},
\end{aligned}
\end{equation}
where $\tilde u$ denotes the solution of \eqref{main} with initial
data $\tilde u_0=\bar u+u_0$.
\end{theo}

In particular, Theorem \ref{nonlin} yields the desired result
of nonlinear stability in the undercompressive or mixed case, effectively
completing the one-dimensional stability analysis initiated in
\cite{ZH,MZ3}.

\begin{rem}\label{rate}
Pointwise bound \eqref{pointwise} yields as a corollary
the sharp $L^p$ decay rate
\begin{equation}
\label{Lp}
|\tilde u(x,t)-\bar u^{\des+\delta(t)}(x)|_{L^p}\le C E_0
(1+t)^{-\frac{1}{2}(1-\frac{1}{p})}, \quad 1\le p\le \infty.
\end{equation}
\end{rem}

The main new difficulty in the analysis beyond those faced in
the strictly parabolic case of \cite{HZ}
is to control higher derivatives in the absence of parabolic smoothing.
We accomplish this by a modification of the fixed-point iteration scheme
introduced in \cite{HZ},
changing to an {\it implicit} iteration scheme in
order to avoid loss of derivatives as discussed in Remark \ref{difference}.

This is a standard device in situations of limited regularity, especially
for quasilinear hyperbolic equations.  However, here the situation is
complicated by the nonlocal character of the
defining integral equations, which appears to prevent the standard
treatment of regularity by energy estimates.
The simple resolution to this problem is that, by appropriate choice
of implicit scheme, we obtain a nonlocal system of integral equations that
admits also a {\it local} description in terms of a symmetric
hyperbolic--parablic system amenable to the same type of energy
estimates used to control regularity in the study of Lax and overcompressive
shocks in \cite{MZ4, Z4, Ra}.

%NOTE (not suff. different from \cite{HR} to put in paper,
%but interesting I think):
%A related but lesser difficulty is the need to integrate Sobolev
%estimates obtained by energy estimates with pointwise bounds
%in a common iteration scheme.
%This issue has been finessed in past pointwise analyses \cite{HR, HRZ}
%%NOTE: I think this is an omission in LZe- they didn't even discuss
%%these issues, just assumed continuity of zeta...
%%for this reason I omit the ref. -KZ
%by splitting the iteration into two steps, obtaining $L^p$ bounds
%in the first then bootstrapping to pointwise bounds in the second.
%This amounts to discarding information, however, replacing spatially
%localized pointwise bounds by a crude supremum estimate,
%an extravagance we cannot afford in the argument at hand.
%In the $L^q\to L^p$ analyses of \cite{MZ4, Z4, Ra} there are no
%pointwise bounds and so this issue does not arise.
%Here, we resolve the issue directly in a simple and general way,
%using a crude weighted energy estimate to establish the boundedness
%and continuity of pointwise estimates needed to carry out
%a standard continuous induction argument; see Section \ref{proof}.
%This small but useful addition to the general framework may
%be helpful in other situations as well.

New physical applications beyond those of \cite{HZ} are to
undercompressive waves in MHD and, with slight
modification following \cite{LRTZ}, to weak detonation waves in
reactive compressible Navier--Stokes equations.
The latter we intend to treat in a future work.

{\bf Plan of the paper.} In Sections \ref{prof} and \ref{lin}, we
recall the basic profile bounds and linearized estimates obtained
in \cite{MZ3, ZH, HRZ}, and in Section \ref{convolution} the
convolution estimates established in \cite{HZ}. In Section
\ref{energy}, we recall (a slight modification of) an auxiliary
energy, or ``hyperbolic--parabolic damping'' estimate, established
in \cite{MZ4, Z4, Ra}, along with a more standard weighted $H^5$
estimate. In Sections \ref{scheme}, \ref{localexistence}, and
\ref{proof}, we carry out the main work of the paper, introducing
a crucial implicit version of the iteration scheme described in
\cite{HZ}, establishing local existence by $H^5$ energy estimate,
and then showing by a combination of estimates like those of
\cite{HZ} and \cite{Z4} that this is contractive in an appropriate
norm encoding the claimed rates of decay.

\section{Profile facts}\label{prof}

We first recall the profile analysis carried out in
\cite{MZ3, HRZ}, a slight generalization of Corollary 1.2, \cite{ZH},
which in turn generalizes results of \cite{MP} in the strictly
parabolic case.
 Profile $\bar u (x)$ satisfies the standing-wave
ordinary differential equation (ODE)
\begin{equation}
B(\bar u) \bar u'=F(\bar u)-F(u_-). \label{ODE}
\end{equation}
Considering the block structure of $B$, this can be written as:
\begin{equation}
F^I(u^I, u^{II})\equiv F^I(u_-^I, u_-^{II})\label{eq1}
\end{equation}
and
\begin{equation}
b_1(u^I)' + b_2(u^{II})'= F^{II}(u^I, u^{II}) - F^{II}(u_-^I,
u_-^{II}). \label{eq2}
\end{equation}

\begin{lem}[\cite{MZ3, HRZ}] \label{profilefact}
Given (H1)--(H3), the endstates $u_\pm$ are hyperbolic rest points
of the ODE determined by (\ref{eq2}) on the $r$-dimensional
manifold (\ref{eq1}), i.e., the coefficients of the linearized
equations about $u_\pm$, written in local coordinates, have no
center subspace. In particular, under regularity (H0),
\begin{equation}\label{expdecay}
D_x^j D_\delta^i(\bar
u^\delta(x)-u_{\pm})=\mathbf{O}(e^{-\alpha|x|}), \quad \alpha>0,
\, 0\le j\le 5, \,i=0,1,
\end{equation}
as $x\rightarrow\pm\infty$.
\end{lem}

%CHANGED: needed update to higher regularity.-K
\begin{proof}
By (H1), \eqref{eq1} may be solved for $u^I=h(u^{II})$,
reducing the problem to an ODE in $u^{II}$.
 Under assumptions (H1)--(H3),
this is a nondegenerate ODE of which
$u^{II}_\pm$ are hyperbolic rest points; see \cite{MZ3, Z4}.
The family $\bar u^\delta$ is thus the intersection of the
unstable manifold at $u_-$ with the stable manifold at $u_+$, both
of which are $C^5$ by (H0) and standard invariant manifold theory.
This intersection is transversal as a consequence of
($\mathcal{D}$), \cite{MZ3}, hence $\bar u^\delta$ is $C^5$ in $\delta$ 
by the Implicit Function Theorem, 
and $C^6$ in $x$ by (H0) and the defining ODE. 
Finally, \eqref{expdecay} follows
from hyperbolicity of $u_\pm$, by standard ODE estimates on
the ODE and its variations about $\bar u^\delta$. For
further details, see \cite{MZ3, Z4}.
\end{proof}
%ENDCHANGED

\section{Linearized estimates}\label{lin}

We next recall some linear theory from \cite{MZ3, ZH}.
Linearizing (\ref{main}) about $\ubar^{\delta_*}(\cdot)$,
$\delta_*$ to be determined later, gives
\begin{equation}
v_t=L^{\delta_*}v:=-(A^{\delta_*}v)_x+(B^{\delta_*}v_x)_x,
\label{linearov}
\end{equation}
with
\begin{equation}
B^{\delta_*}(x):= B(\ubar^{\delta_*}(x)), \quad
A^{\delta_*}(x)v:=
dF(\ubar^{\delta_*}(x))v-dB(\ubar^{\delta_*}(x))v\ubar^{\delta_*}_x.
\label{AandBov}
\end{equation}
Denoting $A^\pm := A(\pm \infty)$,  $B^\pm:= B(\pm \infty)$, and
considering  Lemma \ref{profilefact}, it follows that
\begin{equation}
|A^{\delta_*}(x)-A^-|= \mathbf {O} (e^{-\eta |x|}), \quad
|B^{\delta_*}(x)-B^-|= \mathbf {O} (e^{-\eta |x|})
\label{ABboundsov}
\end{equation}
as $x\to -\infty,$ for some positive $\eta.$ Similarly for $A^+$
and $B^+,$ as $x\to +\infty.$ Also $|A^{\delta_*}(x)-A^\pm|$ and
$|B^{\delta_*}(x)-B^\pm|$ are bounded for all $x$.

 Define the {\it (scalar) characteristic speeds} $a^\pm_1<
\cdots < a_n^\pm$ (as above) to be the eigenvalues of $A^\pm$, and
the left and right {\it (scalar) characteristic modes} $l_j^\pm$,
$r_j^\pm$ to be corresponding left and right eigenvectors,
respectively (i.e., $A^\pm r_j^\pm = a_j^\pm r_j^\pm,$ etc.),
normalized so that $l^+_j \cdot r^+_k=\delta^j_k$ and $l^-_j \cdot
r^-_k=\delta^j_k$. Following Kawashima \cite{Kaw}, define
associated {\it effective scalar diffusion rates}
$\beta^\pm_j:j=1,\cdots,n$ by relation
\begin{equation}
\left(
\begin{matrix}
\beta_1^\pm &&0\\
&\vdots &\\
0&&\beta_n^\pm
\end{matrix}
\right) \quad = \hbox{diag}\ L^\pm B^\pm R^\pm, \label{betaov}
\end{equation}
where $L^\pm:=(l_1^\pm,\dots,l_n^\pm)^t$,
$R^\pm:=(r_1^\pm,\dots,r_n^\pm)$ diagonalize $A^\pm$.
%One of the
%results of the our hypotheses is the fact that
%$$\beta_i^\pm > 0.$$

Assume for  $A$ and $B$ the block structures:
$$A=\left(\begin{matrix}A_{11}\quad A_{12}\\A_{21}\quad A_{22}\end{matrix}\right),
B=\left(\begin{matrix}0& 0\\B_{21}& B_{22}\end{matrix}\right).$$

%%%%%%%%%%%%%%%%%%%%%%%%%%%%%%%%%%%%%%%%%%%%%%%%%%%%%%%%%%%%%%%%%%%%%%%%%%%%%%%%%%%

Also, let $a^{*}_j(x)$, $j=1,\dots,(n-r)$ denote the eigenvalues
of
$$
A_{*}:= A_{11}- A_{12} B_{22}^{-1}B_{21}, \label{A*}
$$
with $l^*_j(x)$, $r^*_j(x)\in \mathbb{R}^{n-r}$ associated left
and right eigenvectors, normalized so that $l^{*t}_jr_j\equiv
\delta^j_k$. More generally, for an $m_j^*$-fold eigenvalue, we
choose $(n-r)\times m_j^* $ blocks $L_j^*$ and $R_j^*$ of
eigenvectors satisfying the {\it dynamical normalization}
$$
L_j^{*t}\partial_x R_j^{*}\equiv 0,
$$
along with the usual static normalization $L^{*t}_jR_j\equiv
\delta^j_kI_{m_j^*}$; as shown in Lemma 4.9, \cite{MZ1}, this
may always be achieved with bounded $L_j^*$, $R_j^*$. Associated
with $L_j^*$, $R_j^*$, define extended, $n\times m_j^*$ blocks
$$
\mathcal{L}_j^*:=\left(\begin{matrix} L_j^* \\
0\end{matrix}\right), \quad \mathcal{R}_j^*:=
\left(\begin{matrix} R_j^*\\
-B_{22}^{-1}B_{21} R_j^*\end{matrix}\right). \label{CalLR}
$$
%
%bridge here: these corr to...
Eigenvalues $a_j^*$ and eigenmodes $\mathcal{L}_j^*$,
$\mathcal{R}_j^*$ correspond, respectively, to short-time
hyperbolic characteristic speeds and modes of propagation for the
reduced, hyperbolic part of degenerate system (\ref{main}).

Define local, $m_j\times m_j$ {\it dissipation coefficients}
$$
\eta_j^*(x):= -L_j^{*t} D_* R_j^* (x), \quad j=1,\dots,J\le n-r,
\label{eta}
$$
where
$$
\aligned &{D_*}(x):=
 \, A_{12}B_{22}^{-1} \Big[A_{21}-A_{22} B_{22}^{-1} B_{21}+ A_{*}
B_{22}^{-1} B_{21} + B_{22}\partial_x (B_{22}^{-1} B_{21})\Big]
\endaligned
\label{D*}
$$
is an effective dissipation analogous to the effective diffusion
predicted by formal, Chapman--Enskog expansion in the (dual)
relaxation case.

The {\it Green distribution} (fundamental solution) associated
with (\ref{linearov}) is defined by
\begin{equation}
G(x,t;y):= e^{L^{\delta_*}t}\delta_y (x). \label{ov3.7}
\end{equation}
or, equivalently,
$$
G_t- L^{\delta_*}G=0, \quad \lim_{t\to 0^+}G(x,t;y)=\delta_y(x).
$$
Recalling the standard notation $ \textrm{errfn} (z) :=
\frac{1}{2\pi} \int_{-\infty}^z e^{-\xi^2} d\xi, $ we have the
following pointwise description.

\begin{prop}[\cite{MZ3}]\label{greenbounds}
Under the assumptions of Theorem \ref{nonlin}, the Green distribution
$G(x,t;y)$ associated with the linearized equations
\eqref{linearov} may be decomposed as $G=H+E+\tilde G$, where,
for $y\le 0$:
\begin{equation}
\begin{aligned} H(x,t;y)&:= \sum_{j=1}^{J} a_j^{*-1}(x) a_j^{*}(y)
\mathcal{R}_j^*(x) \zeta_j^*(y,t) \delta_{x-\bar a_j^* t}(-y)
\mathcal{L}_j^{*t}(y)\\
&= \sum_{j=1}^{J} \mathcal{R}_j^*(x) \mathcal{O}(e^{-\eta_0 t})
\delta_{x-\bar a_j^* t}(-y) \mathcal{L}_j^{*t}(y),
\end{aligned}
\label{multH}
\end{equation}
where the averaged convection rates $\bar a_j^*= \bar a_j^*(x,t)$
in (\ref{multH}) denote the time-averages over $[0,t]$ of
$a_j^*(x)$ along backward characteristic paths $z_j^*=z_j^*(x,t)$
defined by
$$
dz_j^*/dt= a_j^*(z_j^*), \quad z_j^*(t)=x,
$$
the dissipation matrix $\zeta_j^*=\zeta_j^*(x,t)\in
\mathbb{R}^{m_j^*\times m_j^*}$ is defined by the {\it dissipative
flow}
$$
d\zeta_j^*/dt= -\eta_j^*(z_j^*)\zeta_j^*, \quad
\zeta_j^*(0)=I_{m_j},
$$
and $\delta_{x-\bar a_j^* t}$ denotes  Dirac distribution centered
at $x-\bar a_j^* t$.
\begin{equation}\label{E}
E(x,t;y)=\sum_{j=1}^\ell
\frac{\partial \bar u^\delta(x)}{\partial \delta_j}_{|\delta=\d}e_j(y,t),
\end{equation}
\begin{equation}\label{e}
  e_j(y,t)=\sum_{a_k^{-}>0}
  \left(\textrm{errfn }\left(\frac{y+a_k^{-}t}{\sqrt{4\beta_k^{-}t}}\right)
  -\textrm{errfn }\left(\frac{y-a_k^{-}t}{\sqrt{4\beta_k^{-}t}}\right)\right)
  l_{jk}^{-}(y)
\end{equation}
for $y\le 0$ and symmetrically for $y\ge 0$, with
%NOTE: e_j continuuous at zero, though we don't need it....
\begin{equation}\label{ljkbounds}
|l_{jk}^\pm|\le C, \qquad
|(\partial/\partial y)l_{jk}^\pm|\le Ce^{-\eta |y|},
\end{equation}
and
\begin{equation}\label{Gbounds}
\begin{aligned}
|\partial_{x,y}^\alpha &\tilde G(x,t;y)|\le  Ce^{-\eta(|x-y|+t)}\\
& +\quad C(t^{-|\alpha|/2}+
|\alpha_y| e^{-\eta|y|} +|\alpha_x| e^{-\eta|x|})
\Big( \sum_{k=1}^n
t^{-1/2}e^{-(x-y-a_k^{-} t)^2/Mt} e^{-\eta x^+} \\
&+
\sum_{a_k^{-} > 0, \, a_j^{-} < 0}
\chi_{\{ |a_k^{-} t|\ge |y| \}}
t^{-1/2} e^{-(x-a_j^{-}(t-|y/a_k^{-}|))^2/Mt}
e^{-\eta x^+}, \\
&+
\sum_{a_k^{-} > 0, \, a_j^{+}> 0}
\chi_{\{ |a_k^{-} t|\ge |y| \}}
t^{-1/2} e^{-(x-a_j^{+} (t-|y/a_k^{-}|))^2/Mt}
e^{-\eta x^-}\Big), \\
\end{aligned}
\end{equation}
$0\le |\alpha| \le 2$ for $y\le 0$ and symmetrically for $y\ge 0$,
for some $\eta$, $C$, $M>0$, where
$a_j^\pm$ are as in Theorem \ref{nonlin},  $\beta_k^\pm>0$,
$x^\pm$ denotes the positive/negative
part of $x$,  indicator function $\chi_{\{ |a_k^{-}t|\ge |y| \}}$ is
$1$ for $|a_k^{-}t|\ge |y|$ and $0$ otherwise.
Moreover, all estimates are uniform in the supressed parameter $\d$.
\end{prop}
\begin{proof}
This is a restatement of the bounds established in
\cite{MZ3} for pure undercompressive, Lax,
or overcompressive type profiles; the same argument
applies also in the mixed under--overcompressive case.
Also, though it was not explicitly stated, uniformity with respect to
$\d$ is a straightforward consequence of the argument.
\end{proof}

\begin{cor}\label{eboundscor}
Under the assumptions of Theorem \ref{nonlin} and the notation
of Proposition \ref{greenbounds},
\begin{equation}\label{ebounds}
\begin{aligned}
|e_j(y,t)|&\le C\sum_{a_k^->0}
  \left(\textrm{errfn }\left(\frac{y+a_k^{-}t}{\sqrt{4\beta_k^{-}t}}\right)
  -\textrm{errfn }\left(\frac{y-a_k^{-}t}{\sqrt{4\beta_k^{-}t}}\right)\right),\\
|e_j (y,t) &- e_j (y,+\infty)| \le C \textrm{errfn} (\frac{|y|-at}{M\sqrt{t}}),
\quad \text{some }\, a>0 \\
|\partial_t  e_j(y,t)|&\le C t^{-1/2} \sum_{a_k^->0} e^{-|y+a_k^-t|^2/Mt},\\
|\partial_y  e_j(y,t)|&\le C t^{-1/2} \sum_{a_k^->0} e^{-|y+a_k^-t|^2/Mt}\\
&\quad +
Ce^{-\eta|y|}
  \left(\textrm{errfn }\left(\frac{y+a_k^{-}t}{\sqrt{4\beta_k^{-}t}}\right)
  -\textrm{errfn }\left(\frac{y-a_k^{-}t}{\sqrt{4\beta_k^{-}t}}\right)\right),\\
|\partial_y e_j (y,t) &- \partial_y e_j(y,+\infty)|
 \le  C t^{-1/2} \sum_{a_k^->0} e^{-|y+a_k^-t|^2/Mt} \\
|\partial_{yt}  e_j(y,t)|&\le C
(t^{-1}+t^{-1/2}e^{-\eta|y|}) \sum_{a_k^->0} e^{-|y+a_k^-t|^2/Mt}\\
\end{aligned}
\end{equation}
for $y\le 0$, and symmetrically for $y\ge 0$.
\end{cor}

\begin{proof}
Straightforward calculation using \eqref{e} and \eqref{ljkbounds};
see \cite{MZ3}.
\end{proof}

From now on, let
\be\label{evec}
e:=\bp e_1\\ \vdots \\ e_\ell \ep.
\ee

\begin{cor}[\cite{HZ}]\label{Ifact}
Under the assumptions of Theorem \ref{nonlin} and the notation
of Proposition \ref{greenbounds},
\be\label{Ifacteq}
\int^\infty_{-\infty} e(y,+\infty))
(\partial \bar u^{\delta}/\partial \delta)_{|\des}(y)\,dy =I_\ell
\ee
\end{cor}

\begin{proof}
This follows from the standard fact that
$L^\des (\partial \bar u^{\delta}/\partial \delta)_{|\des}=0$,
hence
$$
\int_{-\infty}^{+\infty}G(x,t;y)
(\partial \bar u^{\delta}/\partial \delta)_{|\des}(y)\, dy\equiv
(\partial \bar u^{\delta}/\partial \delta)_{|\des}(x)\,
$$
which,
together with the fact that
$E = (\partial \bar u^{\delta}/\partial \delta)_{|\des}(x) e(y,t))$
represents the only nondecaying
part of $G(x,t;y)$ under stability criterion ($\mathcal{D}$),
yields
$$
(\partial \bar u^{\delta}/\partial \delta)_{|\des}(x)
\int_{-\infty}^{+\infty}e(x,+\infty ;y)
(\partial \bar u^{\delta}/\partial \delta)_{|\des}(y)\, dy=
(\partial \bar u^{\delta}/\partial \delta)_{|\des}(x)\,
$$
in the limit as $t\to +\infty$.
\end{proof}

\begin{prop}[Parameter-dependent bounds]\label{parambds}
Under the assumptions of Theorem \ref{nonlin} and the notation
of Proposition \ref{greenbounds},
\be\label{Hderivbds}
\partial_\des H(x,t;y)=O(H) + \partial_y O(tH),
\ee
\ba\label{ederivbds}
|\partial_\des  e|&\le C|e|\\
|\partial_\des e_t|&\le C | e_t|,\\
|\partial_\des  (e(y,t)-e(y,+\infty))|&\le C
|(e(y,t)-e(y,+\infty))|\\
|\partial_\des  e_y|&\le C| e_y|,\\
|\partial_\des (e_y(y,t)-e_y(y,+\infty)|&\le
C|(e_y(y,t)-e_y(y,+\infty)|\\
|\partial_\des e_{yt}| &\le C| e_{yt}|,\\
\ea
and
\begin{equation}\label{Gderivbds}
\begin{aligned}
|\partial_\des \partial_{x,y}^\alpha &\tilde G(x,t;y)|\le
C|\partial_{x,y}^\alpha &\tilde G(x,t;y)|,
%Ce^{-\eta(|x-y|+t)}\\
%& +\quad C(t^{-|\alpha|/2}+
%|\alpha_y| e^{-\eta|y|} +|\alpha_x| e^{-\eta|x|})
%\Big( \sum_{k=1}^n
%t^{-1/2}e^{-(x-y-a_k^{-} t)^2/Mt} e^{-\eta x^+} \\
%&+
%\sum_{a_k^{-} > 0, \, a_j^{-} < 0}
%\chi_{\{ |a_k^{-} t|\ge |y| \}}
%t^{-1/2} e^{-(x-a_j^{-}(t-|y/a_k^{-}|))^2/Mt}
%e^{-\eta x^+}, \\
%&+
%\sum_{a_k^{-} > 0, \, a_j^{+}> 0}
%\chi_{\{ |a_k^{-} t|\ge |y| \}}
%t^{-1/2} e^{-(x-a_j^{+} (t-|y/a_k^{-}|))^2/Mt}
%e^{-\eta x^-}\Big), \\
\end{aligned}
\end{equation}
$0\le |\alpha| \le 2$ for $y\le 0$ and symmetrically for $y\ge 0$,
for some $C>0$.
\end{prop}

\begin{proof}
These follow by the same argument used to
establish the parameter-dependent bounds of Proposition
3.11, \cite{TZ1}, using the additional fact that neither speed $s$ nor
endstates $u_\pm$ depend on the choice of $\des$ to obtain
better decay estimates on certain terms.

Bounds \eqref{Hderivbds} and \eqref{ederivbds} follow
by direct calculation, together with the observations
(obtained similarly as bounds established in the proof
of Proposition 3.11, \cite{TZ1}, using parameter-dependent
asymptotic ODE bounds) that
$$
\partial_\des a_j^*, \,
\partial_\des \bar a_j^*, \,
\partial_\des \CalR_J^*, \,
\partial_\des \CalL_J^*
=O(1)
$$
and
\be\label{ljkbds}
\partial_\des \partial_{yt}^\alpha l_jk^\pm(y)=
O(e^{-\eta|y|}), \quad 0\le |\alpha|\le 2.
\ee
Bounds \eqref{Gderivbds} follow by the argument of \cite{TZ1},
but using the fact that $\alpha=a_k^\pm, \, \beta_k^\pm$ (since $u_\pm$)
do not depend on $\des$, hence ``bad'' factors $t\partial_\des \alpha=O(t)$
do not appear, but only factors $O(1)$ or better.
\end{proof}

\br\label{Hderivrmk}
\textup{
The additional factor $t$ in the righthand side of \eqref{Hderivbds}
may be absorbed in time-exponential decay of $H$, i.e.,
$tH$ obeys the same decay bounds as $H$,
but with slightly smaller time-exponential decay rate.
}
\er

\medbreak

\section{Convolution lemmas}\label{convolution}
We shall make use of the following technical lemmas proved
in \cite{HZ, HR}.

\begin{lem}[Linear
estimates I]\label{iniconvolutions} Under the assumptions of
Theorem \ref{nonlin},
\begin{equation}\label{iniconeq}
\begin{aligned}
\int_{-\infty}^{+\infty}|\tilde G(x,t;y)|(1+|y|)^{-3/2}\, dy
&\le C(\theta+\psi_1+\psi_2)(x,t),\\
\int_{-\infty}^{+\infty}|e_t(y,t)|(1+|y|)^{-3/2}\, dy
&\le C(1+t)^{-3/2},\\
\int_{-\infty}^{+\infty}|e(y,t)|(1+|y|)^{-3/2}\, dy
&\le C,\\
 \int^{+\infty}_{-\infty} |e(y,t)-e(y,+\infty)| (1+|y|)^{-3/2}\, dy
&\le C(1+t)^{-1/2},\\
\end{aligned}
\end{equation}
for $0\le t\le +\infty$,
some $C>0$, where $\tilde G$ and
$e$ are defined as in Proposition \ref{greenbounds}.
\end{lem}

\begin{lem}[Nonlinear
estimates I]\label{convolutions} Under the assumptions of Theorem
\ref{nonlin},
\begin{equation}\label{coneq}
\begin{aligned}
\int_0^t\int_{-\infty}^{+\infty}|\tilde G_y(x,t-s;y)|\Psi(y,s)\, dy ds
&\le C(\theta+\psi_1+\psi_2)(x,t),\\
%Don't need
%\int_0^t\int_{-\infty}^{+\infty}|\tilde G_x(x,t-s;y)|\Psi(y,s)\, dy ds
%&\le C(\theta+\psi_1+\psi_2)(x,t),\\
\int_0^{t-1}\int_{-\infty}^{+\infty}|\tilde G_{xy}(x,t-s;y)|\Psi(y,s)\, dy ds
&\le C(\theta+\psi_1+\psi_2)(x,t),\\
\int_0^t\int_{-\infty}^{+\infty}|e_{yt}(y,t-s)|\Psi(y,s)\, dy ds
&\le C(1+t)^{-1},\\
\int_t^{+\infty} \int_{-\infty}^{+\infty}|e_y(y,+\infty)|\Psi(y,s)\, dy
&\le C(1+t)^{-1/2},\\
\int_0^t\int_{-\infty}^{+\infty}
|e_y(y,t-s)- e_y(y,+\infty)| \Psi(y,s)\, dyds
&\le C(1+t)^{-1/2},\\
\end{aligned}
\end{equation}
for $0\le t\le +\infty$,
some $C>0$, where $\tilde G$ and
$e$ are defined as in Proposition \ref{greenbounds} and
\begin{equation}\label{source}
\begin{aligned}
\Psi(y,s)&:=
(1+s)^{1/2}s^{-1/2}(\theta + \psi_1+\psi_2)^2(y,s)\\
&\qquad +
(1+s)^{-1} (\theta+\psi_1+\psi_2)(y,s).
\end{aligned}
\end{equation}
\end{lem}

\begin{lem}[Nonlinear estimates II]\label{occonvolutions}
Under the assumptions of Theorem \ref{nonlin},
\begin{equation}\label{occoneq1}
\begin{aligned}
\int_0^t\int_{-\infty}^{+\infty} |\tilde G_y(x,t-s;y)|
\Phi_1(y,s) \, dy ds
&\le C(\theta+\psi_1+\psi_2)(x,t),\\
%Don't need...
%\int_0^t\int_{-\infty}^{+\infty} \tilde G_{x}(x,t-s;y)
%\Phi_1(y,s) \, dy ds
%&\le C(\theta+\psi_1+\psi_2)(x,t),\\
\int_0^{t-1}\int_{-\infty}^{+\infty}|\tilde G_{xy}(x,t-s;y)|\Phi_1(y,s)\, dy ds
&\le C(\theta+\psi_1+\psi_2)(x,t),\\
\int_0^t\int_{-\infty}^{+\infty} |e_{yt}(y,t-s)|
\Phi_1(y,s) \, dy ds
&\le C(1+t)^{-1},\\
\int_0^t\int_{-\infty}^{+\infty} |e_{y}(y,t-s)|
\Phi_1(y,s) \, dy ds
&\le C(1+t)^{-1/2}\\
\end{aligned}
\end{equation}
and
\begin{equation}\label{occoneq2}
\begin{aligned}
\int_0^t\int_{-\infty}^{+\infty} |\tilde G(x,t-s;y)|
\Phi_2(y,s) \, dy ds
&\le C(\theta+\psi_1+\psi_2)(x,t),\\
\int_0^t\int_{-\infty}^{+\infty} \tilde G_x
\Phi_2(y,s) \, dy ds
&\le C(\theta+\psi_1+\psi_2)(x,t),\\
\int_0^t\int_{-\infty}^{+\infty} |e_{t}(y,t-s)|
\Phi_2(y,s) \, dy ds
&\le C(1+t)^{-3/2},\\
\int_0^t\int_{-\infty}^{+\infty} |e(y,t-s)-e(y,+\infty)|
\Phi_2(y,s) \, dy ds
&\le C(1+t)^{-3/2},\\
\end{aligned}
\end{equation}
for $0\le t\le +\infty$,
some $C>0$, where $\tilde G$ and
$e$ are defined as in Proposition \ref{greenbounds} and
\begin{equation}\label{Phi}
\begin{aligned}
\Phi_1(y,s)&:=
e^{-\eta|y|}s^{-1/2}(\theta+\psi_1+\psi_2)(y,s)
\le
Ce^{-\eta|y|/2}s^{-1/2}(1+s)^{-1},\\
%NOTE: It is
%e^{-\eta|y|}(1+s)^{-1/2}s^{-1/2}(1+s)^{1/2}(\theta+\psi_1+\psi_2)(y,s),\\
\Phi_2(y,s)&:=
e^{-\eta|y|}(1+s)^{-3/2}.\\
\end{aligned}
\end{equation}
\end{lem}

%CHANGED: was mislabeled..-K
\begin{lem}[Linear estimates II]\label{Hlinear}
%ENDCHANGED
Under the assumptions of Theorem \ref{nonlin},
if $|v_0(x)|$, $|\partial_x v_0(x)|\le E_0 (1+|x|)^{-\frac32}$,
$E_0>0$,  then, for some $\theta>0$,
\begin{equation}\label{Hinit}
\begin{aligned}
\int_{-\infty}^{+\infty} H(x, t; y) v_0(y) dy&\le C E_0 e^{-\theta
t}(1+|x|)^{-\frac32}\\
 \int_{-\infty}^{+\infty} H_x(x, t; y) v_0(y) dy&\le C E_0
e^{-\theta t}(1+|x|)^{-\frac32},
\end{aligned}
\end{equation}
and so both expressions are dominated by $E_0(\psi_1+\psi_2).$
\end{lem}
\begin{proof}
See \cite{HR, HRZ}.
\end{proof}

%CHANGED
\begin{lem}[Nonlinear estimates III]\label{Hnonlinear}
%ENDCHANGED
%CHANGED: wrong estimate!  This is with decay from your paper,
%not ours (we get only L^2 and thus Sobolev bound of (1+t)^{-1/4} only....
%OLD:
%If  \,$|\Upsilon(y,s)|\le s^{-1/2}(\psi_1 +
%\psi_2 + \theta )(y,s)+ s^{-1/2}e^{-\eta |y|},$
Under the assumptions of Theorem \ref{nonlin},
if  \,$|\Upsilon(y,s)|\le s^{-1/4}(\psi_1 +
\psi_2 + \theta )(y,s)+ s^{-1/2}e^{-\eta |y|},$
then
\ba\label{Hf}
 \Big|\int_0^t \int_{-\infty}^{+\infty}  H (x,t-s;y)  \Upsilon(y,s) dyds\Big|
&\le C ( \psi_1 + \psi_2) (x,t),\\
 \Big|\int_0^t \int_{-\infty}^{+\infty}  (H_x-H_y) (x,t-s;y)  \Upsilon(y,s) dyds\Big|
&\le C ( \psi_1 + \psi_2) (x,t),\\
\int_{t-1}^{t}\int_{-\infty}^{+\infty}|\tilde G_{x}(x,t-s;y)|\Upsilon(y,s)\, dy ds
&\le C(\psi_1+\psi_2)(x,t).\\
\ea
\end{lem}
\begin{proof}
%CHANGED
The proof of the first inequality
%ENDCHANGED
is very much similar to that of the similar estimates proved in
\cite{HR}. Here we prove only the part which contains the
convolution of $H$ against $s^{-1/4}\psi_1$.
A typical term in (\ref{Hf}) coming from $\psi_1$  would be
dominated by a term of the form
\begin{equation*}
\begin{aligned}
C\int_0^t e^{-\eta_0(t-s)}\chi(x-\bar a_j^* (t-s), s)& s^{-\frac14}(1+s)^{-\frac12}\\
\times &(1+|x-\bar a_j^* (t-s)- a_i^-s |)^{-\frac12} ds \\=C\int_0^t e^{-\eta_0(t-s)}\chi(x-\bar a_j^* (t-s), s)& s^{-\frac14}(1+s)^{-\frac12}\\
\times &(1+|x- a_i^-t -(\bar a_j^*-a_i^-) (t-s)|)^{-\frac12}
ds;\end{aligned}\end{equation*}
 now, using $ \frac {1}{1+|a+b|}\le
\frac{1+|b|}{1+|a|} $
 the above  would be smaller than
\begin{equation*}
\begin{aligned} C\int_0^t e^{-\eta_0(t-s)}\chi(x-\bar a_j^* (t-s)&, s) s^{-\frac14}(1+s)^{-\frac12}\\
\times &(1+|x- a_i^-t |)^{-\frac12}(|1+(\bar a_j^*-a_i^-)
(t-s)|)^{\frac12} ds. \end{aligned}\end{equation*}
 Notice that
$(\bar a_j^*-a_i^-)\le C$. Now we use half of $\eta$ to neutralize
$t-s$, and to get
$$C(1+|x- a_i^-t |)^{-\frac12} \int_0^t e^{-\frac{\eta_0}2(t-s)}\chi(x-\bar a_j^* (t-s), s) s^{-\frac14}(1+s)^{-\frac12}ds$$
On the one hand, this is obviously dominated by
$$ C(1+|x- a_i^-t |)^{-\frac12}(1+t)^{-\frac12}.$$
which is absorbed  by a factor of $\psi_1,$ if $x\in [a_1^-t,
a_n^+t].$ On the other hand, $\chi(x-\bar a_j^* (t-s), s)$ in the
integral means that the above expression would vanish if $x>Mt$,
for some fixed $M$. Therefore the above expression is always less
than a factor of $\psi_1+\psi_2.$

The second inequality follows by an identical proof,
using $H_x-H_y \sim H$.
The proof of the third inequality is similar, using the fact that
$$
C|\tilde G_y(t-s)|\le
e^{-\theta (t-s)}
(t-s)^{-1} e^{-|x-y|^2/C(t-s)}
$$
for $t-1\le s\le t$; see
\cite{MZ2} for similar calculations.
\end{proof}

\section{Auxiliary energy estimate}\label{energy}

We shall require the following auxiliary energy
estimate adapted essentially unchanged from \cite{MZ4, Z4, Ra}.
Let
\be\label{viscous}
\tilde u_t + f(\tilde u)_x-(B(\tilde u)\tilde u_x)_x=(\partial \bar u/
\partial \delta)|_{\delta_*}(x) \gamma(t),
\ee
\be\label{pert}
u:=\tilde u -\bar u^{\delta_*+\delta(t)}.
\ee

\begin{lem} [\cite{MZ4, Z4, Ra}] \label{aux}
Under the hypotheses of Theorem \ref{nonlin} let $u_0\in H^5$, and
suppose that, for $0\le t\le T$, the suprema of $ |\dot\delta|$
and $|\gamma|$ and the $W^{3,\infty}$ norm of the solution
$u=(u^I,u^{II})^t$ of \eqref{viscous}, \eqref{pert} each remain
bounded by a sufficiently small constant $\zeta>0$. Then, for all
$0\le t\le T$, \be \label{Ebounds} |u(t)|_{H^5}^2\le C e^{-\theta
t} |u(0)|^2_{H^5} +C\int_0^t e^{-\theta_2 (t-\tau )}(|u|_{L^2}^2+
|\dot\delta|^2 + |\gamma|^2)(\tau)\, d\tau. \ee
\end{lem}

\begin{proof}
This follows exactly as in the $\gamma\equiv 0$ case treated
in \cite{MZ4, Z4, Ra}, observing that term
$$
(\partial \bar u/ \partial \delta)|_{\delta_*}(x) \gamma(t)
$$
is of the same form as terms
$$
(\partial \bar u/ \partial \delta)|_{\delta_*}(x) \cdot\delta(t)
$$
already arising in the nonlinear perturbation equations in the former case.
\end{proof}

We require also the following much cruder estimate adapted from \cite{HR}.

\begin{lem} [\cite{HR}] \label{weightest}
Under the hypotheses of Theorem \ref{nonlin} let
$E_0:=\|(1+|x|^2)^{-3/4}u_0(x)\|_{H^5}<\infty$, and suppose that,
for $0\le t\le T$, the suprema of $ |\dot\delta|$ and $|\gamma|$
and the $W^{3,\infty}$ norm of the solution $u=(u^I,u^{II})^t$ of
\eqref{viscous}, \eqref{pert} each remain bounded by some constant
$C>0$. Then, for all $0\le t\le T$, some $M=M(C)>0$, \ba
\label{weightedEbounds} \|(1+|x|^2)^{-3/4}u(x,t)\|_{H^5}^2&\le
Me^{Mt}\Big( E_0 +\int_0^t (|\dot\delta|^2 + |\gamma|^2)(\tau)\,
d\tau\Big). \ea
\end{lem}

\begin{proof}
This follows by standard Friedrichs symmetrizer estimates carried
out in the weighted $H^5$ norm. (Recall, these plus several more
complicated estimates are used in the proof of Lemma \ref{aux}.)
\end{proof}

\br\label{zetareg} \textup{ An immediate consequence of Lemma
\ref{weightest}, by Sobolev embedding and equation
\eqref{viscous}, is that, if
$E_0:=\|(1+|x|^2)^{-3/4}u_0(x)\|_{H^5}$, $\|u\|_{H^5}$, $|\dot
\delta(\cdot)|$ and $|\gamma(\cdot)|$ are uniformly bounded on
$0\le t\le T$, then
$$
 |(1+|x|^2)^{-3/4}u(x,t)|, \,
|(1+|x|^2)^{-3/4}u_t(x,t)|$$
are uniformly bounded on $0\le t\le T$ as well.
}
\er

\section{Fixed-point iteration scheme}\label{scheme}

We now introduce the fixed-point iteration scheme by which
we shall simultaneously construct and estimate the solution
of the perturbed shock problem.

Our starting point, similarly as in \cite{HZ}, is the observation that
\begin{equation}\label{ocpert}
u(x,t):=\tilde u(x, t)-\bar u^{\delta_*+\delta(t)}(x)
\end{equation}
satisfies
\begin{equation}
\label{mixedperteq}
u_t-L^\des u=Q^\des(u)_x+\dot \delta (t)
(\partial \bar u^\delta/\partial \delta)_{|\des}
+ R^\des(\delta, u)_x +
S^\des(\delta, \delta_t),
\end{equation}
where
%\begin{equation}\label{mixedQ}
%Q^\d =\mathcal{O}( |u|^2),
%\end{equation}
\begin{equation}\label{Q}
\begin{aligned}
Q^\des(u,u_x)&:=
\Big(f(\bar u^{\des})+ A(\bar u^{\des})u
-f(\bar u^{\des}+u) \Big)
+
\Big(B(\bar u^{\des}+u) - B(\bar u^{\des})\Big)u_x
\\
&=\mathcal{O}(|u|^2 + |u||u_x|),\\
Q^\des(u,u_x)_x&=\mathcal{O}(|u||u_x|+|u_x|^2 + |u||u_{xx}|),\\
Q^\des(u,u_x)_{xx}&=\mathcal{O}(|u||u_{xx}|+|u_x|^2 + |u||u_{xxx}|+
|u_x||u_{xx}|),\\
\end{aligned}
\end{equation}
\ba\label{mixedR}
R^\des&:=
\Big(A(\bar u^\des(x))-A(\bar u^{\des+\delta(t)}(x))\Big)u
=\mathcal{O}( e^{-\eta|x|}|\delta||u| ),\\
R^\des_x &= \mathcal{O}( e^{-\eta|x|}|\delta||u_x| ),\\
R^\des_{xx} &= \mathcal{O}( e^{-\eta|x|}|\delta||u_{xx}| ),\\
\ea{equation}
and
\begin{equation}\label{S}
\begin{aligned}
\CalS^\des &:=
\dot \delta \Big( (\partial \bar u^\delta/\partial \delta)_{|\des+\delta(t)}-
(\partial \bar u^\delta/\partial \delta)_{|\des}
\Big)
=\mathcal{O}( e^{-\eta|x|}|\dot \delta|| \delta| )
\end{aligned}
\end{equation}
so long as $|u|$ remains bounded,
by Taylor's Theorem together with \eqref{expdecay}.

Accordingly, for given $\des^{n-1}, \de^{n-1}(\cdot)$,
define $u^n$ to be the solution of
\begin{equation}
\begin{aligned}
u^n(x,t)=&\int_{-\infty}^{+\infty}(H^{n-1}+ \tilde
G^{n-1})(x, t; y) u^{{n-1}}_0(y)dy\\
 &+\int_0^t\int_{-\infty}^{+\infty} (H^{n-1}+
\tilde G^{n-1})(x, t-s; y) \\
&\quad \times S^{\des^{n-1}}(\de^{n-1},\dde^{n-1})(y,s)dyds\\
&-\int_0^t\int_{-\infty}^{+\infty}\tilde G^{n-1}_y(x,
t-s; y)
\\ &\quad \times
\big(Q^{\des^{n-1}}(u^n, u^n_x)+R^{\des^{n-1}}(\de^{n-1},
u^n)\big)(y,s)dyds\\
&+\int_0^t\int_{-\infty}^{+\infty}H^{n-1}(x,
t-s; y)
\\ &\quad \times
\big(Q^{\des^{n-1}}(u^n, u^n_x)+R^{\des^{n-1}}(\de^{n-1},
u^n)\big)_y(y,s)dyds
\end{aligned}\label{un}
\end{equation}
where $H^{n-1}, \tilde G^{n-1}$ are the parts of Green distribution
$G^{n-1}$ of linearized equation around $\bar u^{\des^{n-1}},$ and
\be\label{inidat}
u_0^{n-1}=\tilde u_0-\bar u^{\des^{n-1}}.
\ee
Further, set
\begin{equation}\label{deltan}
 \begin{aligned}
  \delta^n(t)&:=\\
&
-\int^\infty_{-\infty}
(e^{n-1}(y,t)-e^{n-1}(y,+\infty))
u_0^{n-1}(y)\,dy \\
 &
 -\int^t_0\int^{+\infty}_{-\infty} (e^{n-1}(y,t-s)-e^{n-1}(y,+\infty))\\
&\quad \times
S^{\des^{n-1}}(\delta^{n-1},\dot\delta^{n-1})(y,s)\,dy\,ds \\
 &
 +\int^t_0\int^{+\infty}_{-\infty} (e^{n-1}_{y}(y,t-s)-e^{n-1}_y(y,+\infty))\\
&\quad \times
(Q^{\des^{n-1}}(u^n, u^n_x)+ R^{\des^{n-1}}(\delta^{n-1}, u^n))(y,s)\,dy\,ds \\
 &
 +\int_t^{+\infty}\int^{+\infty}_{-\infty} e^{n-1}(y,+\infty)
S^{\des^{n-1}}(\delta^{n-1},\dot\delta^{n-1}) (y,s)\,dy\,ds, \\
 &
 -\int_t^{+\infty}\int^{+\infty}_{-\infty} e^{n-1}_y(y,+\infty)\\
&\quad \times
(Q^{\des^{n-1}}(u^n, u^n_x)+R^{\des^{n-1}}
(\delta^{n-1}, u^n)) (y,s)\,dy\,ds \\
\end{aligned}
\end{equation}
%TODO: check above.
and
\begin{equation}
\begin{aligned}
\des^n:=&\des^{n-1}
+\int_{-\infty}^{+\infty}e^{n-1}(y,+\infty) u^{{n-1}}_0(y)dy\\
 &+\int_0^{+\infty}\int_{-\infty}^{+\infty} e^{n-1}(y, +\infty) S^{\des^{n-1}}(\de^{n-1},\dde^{n-1})(y,s)dyds\\
&-\int_0^{+\infty}\int_{-\infty}^{+\infty}e_y^{n-1}(y,
+\infty)\\
&\quad \times
\big(Q^{\des^{n-1}}(u^n, u^n_x)+R^{\des^{n-1}}(\de^{n-1},
u^n)\big)(y,s)dyds.
\end{aligned}\label{deltasn}
\end{equation}

Define an associated iteration map $\CalT$, formally, by
\be\label{CalT}
(\delta^n, \delta_*^n)=
\CalT (\delta^{n-1}, \delta_*^{n-1}).
\ee

\begin{lem}\label{soln}
Under \eqref{un}--\eqref{deltasn},
$\tilde u^n := u^n +\bar u^{\des^{n-1}+\de^{n-1}}$
satisfies
\begin{equation}\label{unequation}
\begin{aligned}
\tilde u^n_t + f(\tilde u^n)_x - (B(\tilde u^n)\tilde
u^n_x)_x&=\big(\dot\delta^n(t) -\dot\delta^{n-1}(t)\big)\frac{\partial \bar
u^\delta}{\partial \delta}|_{\delta_*^{n-1}},\\
\end{aligned}
\end{equation}
\begin{equation}
\label{dotdeltan}
 \begin{aligned}
  \dot \delta^n (t) &=-\int^\infty_{-\infty}e^{n-1}_t(y,t)
u_0^{n-1}(y)\,dy \\
&\qquad
  - \int^t_0\int^{+\infty}_{-\infty} e^{n-1}_{t}(y,t-s)
%\\&\quad \times
S^{\des^{n-1}}(\delta^{n-1},\dot\delta^{n-1})(y,s)\,dy\,ds \\
 &\qquad
  +\int^t_0\int^{+\infty}_{-\infty} e^{n-1}_{yt}(y,t-s)
\\&\quad \times
(Q^{\des^{n-1}}(u^n, u^n_x)+ R^{\des^{n-1}}(\delta^{n-1}, u^n))(y,s)\,dy\,ds, \\
 \end{aligned}
\end{equation}
with initial data
$\tilde u^n(\cdot, 0)=
\tilde u_0
+
(\bar u^{\des^{n-1}+\de^{n-1}(0)} -\bar u^{\des^{n-1}}),$
and therefore satisfies
\eqref{main}
with initial data $\tilde u_0$ if and only if
\be\label{itcond}
(\delta^n, \delta_*^n)= (\delta^{n-1}, \delta_*^{n-1}),
\ee
i.e., $(\delta^n, \delta_*^n)$ is a fixed point of $\CalT$,
in which case also $\delta(0)=\delta(+\infty)=0$.
\end{lem}

\begin{proof}
Equation \eqref{dotdeltan} follows immeditately upon differentiation of
\eqref{deltan}.
From \eqref{deltan}, we obtain, further, that $\delta^n(+\infty)=0$,
and
\begin{equation}
\begin{aligned}
\de^n(t)-\de^n(0)=&-\int_{-\infty}^{+\infty}e^{n-1}(y,t) u^{{n-1}}_0(y)dy\\
 &-\int_0^t\int_{-\infty}^{+\infty} e^{n-1}(y, t-s) S^{\des^{n-1}}(\de^{n-1},\dde^{n-1})(y,s)dyds\\
&-\int_0^t\int_{-\infty}^{+\infty}e^{n-1}(y,
t-s)\\
&\quad \times
\big(Q^{\des^{n-1}}(u^n, u^n_x)+R^{\des^{n-1}}(\de^{n-1},
u^n)\big)_y(y,s)dyds.
\end{aligned}\label{shifteddeltan}
\end{equation}
Setting $t=+\infty$ in \eqref{shifteddeltan}, and
comparing with \eqref{deltasn}, we find, therefore, that
$\delta^n(0)=\delta^n(+\infty)=0$
if and only if $\delta^n_*=\delta^{n-1}_*$.

From (\ref{un}) and (\ref{deltan}) we conclude that
\begin{equation}
\begin{aligned}
u^n(x,t)=&\int_{-\infty}^{+\infty}
G^{n-1}(x, t; y) u^{{n-1}}_0(y)dy\\
 &+\int_0^t\int_{-\infty}^{+\infty}  G^{n-1}(x, t-s; y) S^{\des^{n-1}}(\de^{n-1},\dde^{n-1})(y,s)dyds\\
&+\int_0^t\int_{-\infty}^{+\infty} G^{n-1}(x, t-s;
y)\\
&\quad \times \big(Q^{\des^{n-1}}(u^n, u^n_x)+R^{\des^{n-1}}(\de^{n-1},
u^n)\big)_y(y,s)dyds\\
&+ (\delta^n(t)-\delta^n(0))\frac{\partial \bar u^\delta}{\partial
\delta}|_{\delta_*^{n-1}}
\end{aligned}
\end{equation}
and thus, by Duhamel's Principle,
\begin{equation}
\begin{aligned}
u^n(t)-L^{\des^{n-1}}u^n&=Q^{\des^{n-1}}(u^n,
u^n_x)_x+R^{\des^{n-1}}(\de^{n-1}, u^n)_x +
S^{\des^{n-1}}(\de^{n-1},\dde^{n-1})  \\
&\quad +\dot
\delta^n(t)\frac{\partial \bar u^\delta}{\partial
\delta}|_{\delta_*^{n-1}}.
\end{aligned}
\end{equation}

Setting $\tilde u^n = u^n +\bar u^{\des^{n-1}+\de^{n-1}}$,
we then obtain \eqref{unequation}
by a straightforward calculation comparing with \eqref{mixedperteq},
with the claimed initial data
\begin{equation}\label{init}
\tilde u^n(\cdot, 0)=
\bar u^{\des^{n-1}+\de^{n-1}(0)}+u^{n-1}_0
=
\tilde u_0
+
(\bar u^{\des^{n-1}+\de^{n-1}(0)} -\bar u^{\des^{n-1}}).
\end{equation}
Note that the righthand side is equal to $\tilde u_0$
if and only if $\delta^{n-1}(0)=0$, or, in case $\delta^n\equiv \delta^{n-1}$
(a fixed point),
if $\delta^n(0)=0$, or equivalently $\delta^n_*=\delta^{n-1}_*$.
\end{proof}

\br\label{difference}
\textup{
Other than a slight notational change $\delta\to \des+\delta$ made
to simplify the exposition, the difference between this iteration
scheme and the one used in \cite{HZ} in the strictly parabolic case
is that we have made it implicit in $u^n$, i.e., $u$ appears everywhere
on the righthand side of the integral equations with index $n$ rather
than $n-1$.
By this change we preserve regularity properties,
as encoded in the nonlinear structure of
equation \eqref{unequation}; see
Lemma \ref{aux}.  By contrast, the explicit version of
\cite{HZ} is not associated with a (favorable) nonlinear equation,
and so would lose derivatives, preventing the iteration scheme from closing.
}
\er

\br\label{inirmk}
\textup{
Note that \eqref{unequation}--\eqref{dotdeltan} form a closed
system for $(u^n, \dot \delta^n)$, in the form of a true {\it Cauchy
problem}; that is, the values of $(u^n, \dot \delta^n)$ at time $T$
depend only on values for $0\le t\le T$, and not on future times.
By \eqref{un}, we have, evidently,
\be\label{ineq}
u^n(\cdot, 0)=u^{n-1}_0.
\ee
}
\er

\section{Local existence}\label{localexistence}

\begin{lem}[$H^5$ local theory] \label{local}
Under the hypotheses of Theorem \ref{nonlin}
let
$$
E_1:=\|u_0(x)\|_{H^5} + \|\delta^{n-1}\|_{B_1} + |\des^{n-1}|
<\infty.
$$
Then, for $T=T(E_1)>0$ sufficiently small
and $C=C(E_1,T)>0$ sufficiently large, there exists
on $0\le t\le T$ a unique solution
$$(u^n,\dot \delta^n)\in L^\infty(H^5(x);t) \times C^{0}(t)
$$
%of \eqref{viscous}, \eqref{pert}, satisfying
of \eqref{unequation}--\eqref{dotdeltan}, satisfying
\be\label{localstab}
%\|u\|_{H^5}(t),\,  |\delta(t)|, \, |\dot \delta(t)|\le CE_1.
\|u^n\|_{H^5}(t),\,  |\dot \delta^n(t)|\le CE_1. \ee
\end{lem}

\begin{proof}
Short-time existence, uniqueness, and stability follow
by (unweighted) energy estimates in $u^n$ similar
to \eqref{weightedEbounds} combined with more straightforward
estimates on $\dot \delta^n$ carried out directly from
integral equation \eqref{dotdeltan},
using a standard (bounded high norm,
contractive low norm) contraction
mapping argument like those described in \cite{Z4, Z5}.
We omit the details.
\end{proof}

\br\label{cauchy}
\textup{
A crucial point is that equations \eqref{unequation}--\eqref{dotdeltan}
depend only on values of $(u^n, \dot\delta^n)$ on the range $t\in [0,T]$;
see Remark \ref{inirmk}.
}
\er

\section{Proof of the Main Theorem}\label{proof}
We are now ready to prove the main theorem.
Define norms
\begin{equation}
\begin{aligned}
|h|_{B_1}&:=
%|f(\theta+\psi_1+\psi_2)^{-1}|_{L^\infty(x,t)} +
|h(t)(1+t)^{1/2}|_{L^\infty(t)}
+
|\dot h(t)(1+t)|_{L^\infty(t)},
\end{aligned}
\label{Beta1norm}
\end{equation}
and
\begin{equation}
\begin{aligned}
|g|_{B_2}&:=
|f(\theta+\psi_1+\psi_2)^{-1}|_{W^{1,\infty}(x,t)}
\end{aligned}
\label{Beta2norm}
\end{equation}
and Banach spaces
\begin{equation}
B_1:=\{h:\, |h|_{B_1}<+\infty\},
\qquad
B_2:=\{g:\, |g|_{B_2}<+\infty\}.
\label{B}
\end{equation}

%NOTE: (Key new part):
\begin{lem}\label{continuation}
Under the hypotheses of Theorem \ref{nonlin}, let
$$(u^n, \dot\delta^n)\in L^\infty(t,H^5(x))\times L^{\infty}(t)$$
satisfy \eqref{unequation}--\eqref{dotdeltan} on $0\le t\le T$,
and define \ba\label{zeta} \zeta(t)&:=\sup_{x, \, 0\le s\le t}
\Big((|u^n|+|u^n_x|)(\theta+\psi_1+\psi_2)^{-1}(x,s)
+|\dot\delta^n(s)(1+s)| \Big). \ea If $\zeta(T)$,
$\|u^{n-1}_0\|_{H^5}$, and
%$\int_0^{+\infty}|\dot \delta^{n-1}(s)|^2\, ds$
$|\delta^{n-1}|_{B_1}$
are bounded by $\zeta_0>0$ sufficiently small, then,
for some $\eps>0$,
(i) the solution $(u^n, \dot \delta^n)$, and thus $\zeta$
extends to $0\le t\le T+\eps$,
and (ii) $\zeta$ is bounded and continuous on $0\le t\le T+\eps$.
\end{lem}

\begin{proof}
By \eqref{ineq} and Lemma \ref{aux},
smallness of $\zeta(T)$,
$|\delta^{n-1}(s)|_{B_1}$, and
$$
\|u^{n-1}_0\|_{H^5}= \|u^{n}(\cdot, 0)\|_{H^5}
%\qquad \hbox{\rm (recall \eqref{ineq})}
$$
(recall \eqref{ineq}) together imply boundedness (and smallness,
though we don't need this) of $\|u\|_{H^5}$ and $|\dot
\delta^n|_{L^\infty}$ on $0\le t\le T$. By Lemma \ref{local}, this
implies existence, boundedness of $\|u\|_{H^5}$,
$\|\delta\|_{W^{1,\infty}}$ on $0\le t\le T+\epsilon$ for
$\eps>0$, and thus, by Remark \ref{zetareg}, boundedness and
continuity of $\zeta$ on $0\le t\le T+\eps$.
\end{proof}

\begin{lem}\label{uexistence}
For $M>0$ and $C_1\ge C>>M$ sufficiently large,
for
$$E_0:=\|(1+|x|^2)^{-3/4}(\tilde u_0- \bar u)\|_{H^5}$$
sufficiently small,  and $|\delta^{n-1}|_{B_1}+ M|\des^{n-1}|\le
2CE_0$, there exist solutions $(u^n, \delta^n, \des^n)$ of
\eqref{un}--\eqref{deltasn} for all $t\ge 0$, satisfying
\be\label{H4unests} |u^n|_{H^5} \le C_1E_0 \ee and
\be\label{unest} |u^n|_{B_2}+
%NOTE: not in induction loop....
|\delta^{n}|_{B_1}+ M|\des^n|\le 2CE_0.
\ee
\end{lem}

\begin{proof}
Define $\zeta$ as in \eqref{zeta}.
Then,
%by Lemma \ref{continuation},
it is sufficient to show that \be\label{claim} \zeta(t)\le CE_0 +
C_*(E_0+\zeta(t))^2 \ee for fixed $C$, $C_*>0$, so long as the
solution $(u^n, \dot\delta^n)\in L^\infty(t,H^5(x))\times
C^{0}(t)$ of \eqref{unequation}--\eqref{dotdeltan} exists and
\be\label{finalbd} \zeta(t)\le (3/2)CE_0, \ee
%$\zeta$ remains sufficiently small,
%in order to obtain by continuous induction
%\be\label{finalbd}
%\zeta(t)\le (3/2)CE_0
%\ee
in order to conclude that solution $(u^n, \dot \delta^n)$ exists
and satisfies \eqref{finalbd} for all $t\ge 0$, provided
$$
E_0< 25/2CC_*
$$
is sufficiently small.

For, by \eqref{ineq} and \eqref{inidat},
and \eqref{expdecay},
\ba\label{inismall}
\|u^{n}(\cdot, 0)\|_{H^5}&= \|u^{n-1}_0\|_{H^5}\\
&= \| \tilde u_0-\bar u^{\des^{n-1}}\|_{H^5}\\
&\le
 \| \tilde u_0-\bar u\|_{H^5}
 + \| \bar u^{\des^{n-1}}- \bar u\|_{H^5}\\
&\le E_0 + c_1|\des^{n-1}|\\
&\le E_0(1 + 2c_1C/M) \\
%&\le (C/2)E_0
\ea
is small, for $E_0$ sufficiently small.
Letting $T$ be the maximum time up to which a solution $(u^n, \dot \delta^n)$
exists and $\zeta \le \zeta_0$ sufficiently small (note: $T\ge 0$ by
the weighted version of \eqref{inismall}, together with \eqref{dotdeltan}),
%NOTE: specifically, easy calc. using bound on $u^n$ in \eqref{dotdeltan}.
we find by Lemma \ref{continuation},  therefore, and
the assumed bounds on $\delta^{n-1}$ and $\dot \delta^{n-1}$,
that $(u^n, \dot \delta^n)$ exists up to $T+\eps$, $\eps>0$,
and that $\zeta$ remains bounded and continuous up to $T+\eps$ as well.
Observing that \eqref{claim} together with
$E_0< 2/(9C^2+6C)$
implies that $\zeta(t)< (3/2)CE_0$ whenever $\zeta(t) \le (3/2)CE_0$,
we find by continuity that  $\zeta(t) \le (3/2)CE_0$ up to $t=T+\eps$
as claimed.
%TODO: expand a bit?
%NOTE: continuity keeps $\zeta$ small as required for continutation step.

By the definition of $\zeta$, we obtain therefore
\be\label{fbd}
|u^n|_{B_2}+ |\dot\delta^n(t)(1+t)|_{L^\infty}\le (3/2)CE_0.
\ee
Thus, it is sufficient to establish first \eqref{claim}, then
afterward, assuming \eqref{finalbd},
\be\label{penbd}
|\delta^n(t)|\le (CE_0/4) (1+t)^{-1/2}
\ee
and
\be\label{lastbd}
|\des^n|\le (CE_0/4M),
\ee
from which we obtain \eqref{unest} by summation with \eqref{fbd};
noting that \eqref{penbd} and \eqref{lastbd}
include the information that
the integral equations for $\delta^n$ and $\des^n$ converge, we obtain also,
by Lemma \ref{soln} and the fact
that $(u^n,\dot \delta^n)$ satisfies \eqref{unequation}--\eqref{dotdeltan}
for all $t\ge 0$, that $(u^n, \delta^n, \des^n)$
satisfies \eqref{un}--\eqref{deltasn} as claimed.
Finally, recalling \eqref{unequation} and applying Lemma \ref{aux}
with $\gamma:=\dot \delta^{n}- \dot \delta^{n-1}$, we obtain
\eqref{H4unests} so long as \eqref{claim} remains valid,
controlling $\|u^n\|_{H^5}$ by integrating the righthand side of
\eqref{Ebounds} and using \eqref{finalbd}, the definition of
$\zeta$, and the assumed bounds on $\dot \delta^{n-1}$.  (We carry
out this last calculation in detail in the following paragraph, in
the course of proving \eqref{claim}).

We now establish \eqref{claim} assuming \eqref{finalbd}. By Lemma
\ref{aux}, and the one-dimensional Sobolev bound
$|u^n|_{W^{3,\infty}}\le c|u^n|_{H^5}$, we have \ba \label{2calc}
|u^n(t)|_{H^5}^2&\le c |u^n(0)|^2_{H^5}e^{-\theta t}\\&\quad
+c\int_0^t e^{-\theta_2 (t-\tau )}(|u^n|_{L^2}^2+ |\dot\delta^n|^2
+|\dot \delta^n- \dot \delta^{n-1}|^2)(\tau)\, d\tau \\
&\le c |u^n(0)|^2_{H^5}e^{-\theta t}\\&\quad +c\int_0^t
e^{-\theta_2 (t-\tau )}(|u^n|_{L^2}^2+ \max\{
|\dot\delta^n|^2, | \dot \delta^{n-1}|^2 \})(\tau)\, d\tau \\
&\le c_2\big(|u^n(0)|^2_{H^5}+ \zeta(t)^2\big) (1+t)^{-1/2}\\
&\le c_2\big(E_0^2(1 + 2c_1C/M)^2
+ (3CE_0/2)^2\big) (1+t)^{-1/2}\\
& \le (C_1E_0)^2 (1+t)^{-1/2},
%changedR \le (C_1E_0)^2,
\ea
for $C_1>0$ sufficiently large $E_0$ sufficiently small,
by \eqref{inismall}, \eqref{finalbd}, and the definition of $\zeta$.
This verifies \eqref{H4unests}, assuming \eqref{finalbd}.

With \eqref{mixedperteq}, \eqref{H4unests} and the resulting
Sobolev estimate $\|u^n\|_{W^{3,\infty}}\le cC_1E_0$,
assumption $|\delta|_{B_1}\le 2CE_0$,
and the definitions of $\zeta$ and $|\cdot|_{B_1}$,
we obtain readily
\ba\label{sourcebds}
|Q^\des +R^\des|&\le c(\zeta^2+4C^2E_0^2)(\Psi+ \Phi_1), \\
|Q^\des_y +R^\des_y|, \,
|Q^\des_{yy} +R^\des_{yy}|&\le
c(\zeta^2+4C^2E_0^2) \Upsilon
\ea
and
\ba\label{derivsourcebds}
|S^\des|,\,
|S^\des_y|&\le
c(\zeta^2+4C^2E_0^2) \Phi_2,
\ea
where $\Phi$, $\Psi$, and $\Upsilon$ are as defined in
Lemmas \ref{iniconvolutions}--\ref{Hnonlinear}.

Expressing $u^n_x$ using \eqref{un} as
\ba\label{unx}
u^n_x(x,t)=&\int_{-\infty}^{+\infty}(H^{n-1}_x+ \tilde
G^{n-1}_x)(x, t; y) u^{{n-1}}_0(y)dy\\
 &+\int_0^t\int_{-\infty}^{+\infty} (H^{n-1}_x-H^{n-1}_y +
\tilde G^{n-1}_x)(x, t-s; y) \\
&\quad \times S^{\des^{n-1}}(\de^{n-1},\dde^{n-1})(y,s)dyds\\
 &-\int_0^t\int_{-\infty}^{+\infty} H^{n-1} (x, t-s; y)
%\\ &\quad \times
S^{\des^{n-1}}_y(\de^{n-1},\dde^{n-1})(y,s)dyds\\
&-\int_0^{t-1}\int_{-\infty}^{+\infty}\tilde G^{n-1}_{yx}(x,
t-s; y)
\\ &\quad \times
\big(Q^{\des^{n-1}}(u^n, u^n_x)+R^{\des^{n-1}}(\de^{n-1},
u^n)\big)(y,s)dyds\\
&+\int_{t-1}^{t}\int_{-\infty}^{+\infty}\tilde G^{n-1}_{x} (x, t-s; y)
\\ &\quad \times
\big(Q^{\des^{n-1}}(u^n, u^n_x)+R^{\des^{n-1}}(\de^{n-1},
u^n)\big)_y(y,s)dyds\\
&+\int_0^t\int_{-\infty}^{+\infty}(H^{n-1}_x-H^{n-1}_y)(x,
t-s; y)
\\ &\quad \times
\big(Q^{\des^{n-1}}(u^n, u^n_x)+R^{\des^{n-1}}(\de^{n-1},
u^n)\big)_y(y,s)dyds,\\
&-\int_0^t\int_{-\infty}^{+\infty}H^{n-1}(x,
t-s; y)
\\ &\quad \times
\big(Q^{\des^{n-1}}(u^n, u^n_x)+R^{\des^{n-1}}(\de^{n-1},
u^n)\big)_{yy}(y,s)dyds,\\
\ea
and applying Lemmas \ref{iniconvolutions}--\ref{Hnonlinear}
to \eqref{un}, \eqref{unx}, and \eqref{dotdeltan},
we thus obtain \eqref{claim} as claimed.

Likewise, we obtain easily \eqref{penbd} from
\eqref{deltan} and \eqref{inismall}, using
Lemmas \ref{iniconvolutions}--\ref{Hnonlinear}
and the definitions of $\zeta$ and $|\cdot|_{B_1}$.

Thus, it remains only to establish \eqref{lastbd}.
This is more delicate, due to the appearance of $M$ in the
denominator of the righthand side, and depends on the
key fact that estimate $\des^n$ of the asymptotic shock location
is to linear order insensitive to the initial guess $\delta^{n-1}$.
To see this, decompose the expression \eqref{deltasn} for
$\des^n$ into linear and nonlinear parts
\ba\label{deslin}
I&:=\des^{n-1}-\int_{-\infty}^{+\infty}e^{n-1}(y,+\infty) u^{{n-1}}_0(y)dy\\
&=
\Big(\des^{n-1}-\int_{-\infty}^{+\infty}e|_{\des=0}(y,+\infty)
(\bar u -\bar u^{\des^{n-1}})(y)dy \Big)\\
&-\int_{-\infty}^{+\infty}e|_{\des=0}(y,+\infty)
(\tilde u_0 -\bar u)(y)dy \\
&\quad - \int_{-\infty}^{+\infty}(e^{n-1}-e|_{\des=0})(y,+\infty)
u^{n-1}_0(y)dy\\
&=: I_a+I_b+I_c
\ea
and
\ba\label{desnonlin}
 II&:=-\int_0^{+\infty}\int_{-\infty}^{+\infty} e^{n-1}(y, +\infty) S^{\des^{n-1}}(\de^{n-1},\dde^{n-1})(y,s)dyds\\
&-\int_0^{+\infty}\int_{-\infty}^{+\infty}e^{n-1}(y,
+\infty)\\
&\quad \times
\big(Q^{\des^{n-1}}(u^n, u^n_x)+R^{\des^{n-1}}(\de^{n-1},
u^n)\big)_y(y,s)dyds,
\ea
respectively.

%Term $I_c$ vanishes, since $e$ depends only on the endstates
%$u_\pm$, which are independent of $\des$.
By estimates like the previous ones, we readily obtain
$$
|II|\le 2c(2CE_0)^2,
$$
which is $<< CE_0/4M$ for $E_0$ sufficiently small.
Likewise, $|I_c|\le c|\des|E_0$, by \eqref{inismall},
\eqref{ederivbds}, and the Mean Value Theorem,
hence is $<<CE_0/4M$ for $E_0$ sufficiently small
(recall that we assume $|\des|\le 2CE_0$),
and
\ba
|I_b|&\le c\|\tilde u_0-\bar u\|_{L^1}\\
&\le  c_2 \|(1+|x|^2)^{-3/4}(\tilde u_0-\bar u)\|_{H^4}\\
&\le
c_2E_0,
\ea
hence is $<<CE_0/4M$ for $C>0$ sufficiently large.

Finally, Taylor expanding, and recalling \eqref{expdecay} and
\eqref{Ifacteq}, we obtain \ba \label{deslina}I_a&=
\des^{n-1}-\des^{n-1}\int_{-\infty}^{+\infty}e|_{\des=0}(y,+\infty)
(\partial \bar u^{\des}/\partial \des)|_{\des=0}(y)dy\\
&\quad  + O(|\des|^2)\\
&=O(|\des|^2),\\
\ea
which is also $<<CE_0/4M$ for $E_0$ sufficiently small
(recall that we assume $|\des|\le 2CE_0$).
%NOTE: We establish \eqref{claim}  (merge from \cite{Z4}, \cite{HZ}),
%in the course of which we find \eqref{lastbd} by essentially
%\eqref{diag} below.
Summing, we obtain \eqref{lastbd} for $E_0$ sufficiently small
and $C>0$ sufficiently large, as claimed.
This completes the proof.
\end{proof}

\begin{proof}[Proof of Theorem \ref{nonlin}]
Define now
\begin{equation}\label{scalednorm}
|(i,j)|_*:= |i|_{B_1}+ M|j|,
\quad
(i,j)\in B_1\times \mathbb{R}.
\end{equation}
By Lemma \ref{uexistence}, for $M>0$ sufficiently large, and
$$
E_0:=\|(1+|x|^2)^{-3/4}(\tilde u_0-\bar u)|_{H^5}
$$
and $r>0$ sufficiently small,
$\CalT=(\mathcal{T}_{\delta} , \mathcal{T}_{\des})$
is a well-defined mapping from
$$
B(0,r)\subset B_1\times \mathbb{R}\to
B_1\times \mathbb{R}.
$$
To establish the theorem, therefore, it suffices to establish that
$\CalT$ is a contraction on $B(0,r)$ in the norm $|\cdot|_*$.
For, then, applying Contraction Mapping Theorem, we find that
$$
(\CalT_\delta, \CalT_\des)(\delta, \des)= (\delta, \des)
$$
has a unique solution
$(\delta^n, \des^n) \in B(0,r)\subset B_1\times \mathbb{R}$,
for which the associated $(u^n, \delta^n, \des^n)$ by Lemma \ref{soln} satisfy
$u^n=\tilde u- \bar u^{\des^n+\delta^n(t)}$ with $\tilde u$ a solution
of \eqref{main} with initial data $\tilde u_0$,
and the stated decay estimates follow by \eqref{H4unests} and \eqref{unest}.

That $\mathcal{T}$ is a contraction follows, provided we can
establish on $B(0,r)$ the Lipshitz bounds
\ba\label{lipbds} |\mathcal{T}(\delta,
\des)-\mathcal{T}(\hat\delta, \hat\des)|_* \le \alpha |(\delta,
\des)-(\hat\delta, \hat\des)|_*
 \ea
for some $\alpha < 1.$

Letting $(u^n,\de^n, \des^n)$ satisfy \eqref{un}--\eqref{deltasn}
for $\delta^{n-1}$, $\des^{n-1}$, and
$(\hat u^n,\hat \de^n, \hat \delta_*^n)$ satisfy \eqref{un}--\eqref{deltasn}
with $\delta^{n-1}$, $\des^{n-1}$
replaced by $\hat \delta^{n-1}$, $\hat \delta_*^{n-1}$, define
variations
\be\label{vardef}
\Delta u^n:=\hat u^n-u^n,
\quad
\Delta \delta^n:=\hat \delta^n-\delta^n,
\quad
\Delta \delta_*^n:=\hat \delta^n_*-\delta^n_*
\ee
and
\be\label{varindef}
\Delta \delta^{n-1}:=\hat \delta^{n-1}-\delta^{n-1},
\quad
\Delta \delta_*^{n-1}:=\hat \delta^{n-1}_*-\delta^{n-1}_*.
\ee
Likewise, define $\Delta \tilde G^{n-1}$, $\Delta H^{n-1}$,
$\Delta e^{n-1}$ in the obvious way.

{\it Differential variational equation}.
From \eqref{unequation}, we find after a brief calculation that
$\Delta \tilde u^n$ defined by
$$
\Delta u^n=\Delta \tilde u^n-
\big(\bar u^{\hat \delta_*^{n-1}+\hat \delta^{n-1}(t)}
-\bar u^{\delta_*^{n-1}+ \delta^{n-1}(t)}\big)
$$
satisfies the variational equations associated with
\begin{equation}\label{varunequation}
\begin{aligned}
\tilde u^n_t + f(\tilde u^n)_x - (B(\tilde u^n)\tilde
u^n_x)_x&=\big(\dot\delta^n(t) -\dot\delta^{n-1}(t)\big)\frac{\partial \bar
u^\delta}{\partial \delta}|_{\delta_*^{n-1}},\\
\end{aligned}
\end{equation}
from which we obtain by an energy estimate similar to that of Lemma
\ref{aux} and the observation
$$
|\Delta u^n(0)|^2_{H^4}= |\bar u^{\delta_*^{n-1}}-\bar u^{\hat
\delta_*^{n-1}}|^2_{H^4} \le C|\Delta \delta_*^{n-1}|^2
$$
the bound
\be \label{varEbounds}
\begin{aligned}
|\Delta u^n(t)|_{H^4}^2&\le
C e^{-\theta t} |\Delta u^n(0)|^2_{H^4}\\
&\quad +C\int_0^t e^{-\theta_2 (t-\tau )}(|\Delta u^n|_{L^2}^2+
\max\{|\Delta \dot \delta^{n}|^2, |\Delta \dot \delta^{n-1}|^2\}\\
&\quad
+
|\Delta \delta_*^{n-1}|^2\max\{|\dot \delta^{n}|^2, |\dot \delta^{n-1}|^2\}
)(\tau)\, d\tau\\
&\le
C\int_0^t e^{-\theta_2 (t-\tau )}(|\Delta u^n|_{L^2}^2+
|\Delta \dot \delta^{n}|^2)(\tau)\, d\tau\\
&\quad +C|(\Delta \delta^{n-1}, \Delta
\delta_*^{n-1})|_{*}^2(1+t)^{-1/2},
%changedR C|(\Delta \delta^{n-1}, \Delta \delta_*^{n-1})|_{*}^2,
\end{aligned}
\ee provided $r$ (and so $\sup\|\hat u^n\|_{H^5}$ and
$\sup\|u^n\|_{H^5}$) is sufficiently small, so long as $\|\Delta
u^n\|_{H^4}$ remains sufficiently small. We omit the (standard)
details.

{\it Integral variational equations}.
Applying the quadratic Leibnitz formula $\Delta (fg)=
f\Delta g + \Delta f g$, we obtain
\begin{equation}
\begin{aligned}
\Delta u^n(x,t) =&
\int_{-\infty}^{+\infty}(\Delta H^{n-1}+ \Delta \tilde
G^{n-1})(x, t; y) u^{{n-1}}_0(y)dy\\
&+\int_{-\infty}^{+\infty}(H^{n-1}+ \tilde
G^{n-1})(x, t; y) \Delta u^{{n-1}}_0(y)dy\\
 &+\int_0^t\int_{-\infty}^{+\infty} (\Delta H^{n-1}+
\Delta \tilde G^{n-1})(x, t-s; y) \\
&\quad \times S^{\des^{n-1}}(\de^{n-1},\dde^{n-1})(y,s)dyds\\
 &+\int_0^t\int_{-\infty}^{+\infty} (H^{n-1}+
\tilde G^{n-1})(x, t-s; y) \\
&\quad \times \Delta S(y,s)dyds\\
&-\int_0^t\int_{-\infty}^{+\infty}\tilde \Delta G^{n-1}_y(x,
t-s; y)
\\ &\quad \times
\big(Q^{\des^{n-1}}(u^n, u^n_x)+R^{\des^{n-1}}(\de^{n-1},
u^n)\big)(y,s)dyds\\
&-\int_0^t\int_{-\infty}^{+\infty}\tilde G^{n-1}_y(x,
t-s; y)
\\ &\quad \times
\big(\Delta Q+ \Delta R \big)(y,s)dyds\\
&+\int_0^t\int_{-\infty}^{+\infty}\Delta H^{n-1}(x,
t-s; y)
\\ &\quad \times
\big(Q^{\des^{n-1}}(u^n, u^n_x)+R^{\des^{n-1}}(\de^{n-1},
u^n)\big)_y(y,s)dyds\\
&+\int_0^t\int_{-\infty}^{+\infty}H^{n-1}(x,
t-s; y)
\\ &\quad \times
\big(\Delta Q+\Delta R)\big)_y(y,s)dyds,
\end{aligned}\label{varun}
\end{equation}
\begin{equation}
\label{vardotdeltan}
 \begin{aligned}
  \Delta \dot \delta^n (t) &=-\int^\infty_{-\infty}\Delta e^{n-1}_t(y,t)
u_0^{n-1}(y)\,dy \\
&-\int^\infty_{-\infty}e^{n-1}_t(y,t)
\Delta u_0^{n-1}(y)\,dy \\
&\qquad
  - \int^t_0\int^{+\infty}_{-\infty} \Delta e^{n-1}_{t}(y,t-s)
S^{\des^{n-1}}(\delta^{n-1},\dot\delta^{n-1})(y,s)\,dy\,ds \\
&\qquad
  - \int^t_0\int^{+\infty}_{-\infty} e^{n-1}_{t}(y,t-s)
\Delta S(y,s)\,dy\,ds \\
 &\qquad
  +\int^t_0\int^{+\infty}_{-\infty} \Delta e^{n-1}_{yt}(y,t-s)
\\&\quad \times
(Q^{\des^{n-1}}(u^n, u^n_x)+ R^{\des^{n-1}}(\delta^{n-1}, u^n))(y,s)\,dy\,ds, \\
 &\qquad
  +\int^t_0\int^{+\infty}_{-\infty} e^{n-1}_{yt}(y,t-s)
\\&\quad \times
(\Delta Q+ \Delta R(y,s)\,dy\,ds, \\
 \end{aligned}
\end{equation}
\begin{equation}
\begin{aligned}
\Delta \des^n:=&\Delta \des^{n-1}
+\int_{-\infty}^{+\infty}\Delta e^{n-1}(y,+\infty) u^{{n-1}}_0(y)dy\\
&+\int_{-\infty}^{+\infty}e^{n-1}(y,+\infty) \Delta u^{{n-1}}_0(y)dy\\
 &+\int_0^{+\infty}\int_{-\infty}^{+\infty} \Delta e^{n-1}(y, +\infty) S^{\des^{n-1}}(\de^{n-1},\dde^{n-1})(y,s)dyds\\
 &+\int_0^{+\infty}\int_{-\infty}^{+\infty} e^{n-1}(y, +\infty) \Delta S (y,s)dyds\\
&-\int_0^{+\infty}\int_{-\infty}^{+\infty}\Delta e_y^{n-1}(y, +\infty)\\
&\quad \times
\big(Q^{\des^{n-1}}(u^n, u^n_x)+R^{\des^{n-1}}(\de^{n-1},
u^n)\big)(y,s)dyds\\
&-\int_0^{+\infty}\int_{-\infty}^{+\infty}e_y^{n-1}(y,
+\infty)\\
&\quad \times
\big(\Delta Q + \Delta R \big)(y,s)dyds,
\end{aligned}\label{vardeltasn}
\end{equation}
and similarly for $\Delta \delta^n$, where
\ba\label{varin}
\Delta u_0^{n-1}&=
\bar u^{\delta_*^{n-1}}
-\bar u^{\hat \delta_*^{n-1}}\\
&=\Big(\frac{\partial \bar u^\delta_*}{\partial \delta_*}\Big)|_{\delta_*=\delta_*^{n-1}}\Delta \delta_*^{n-1}+
O(|\Delta \delta_*^{n-1}|^2e^{-\eta|x|})
\\
& = O(|\Delta \delta_*^{n-1}e^{-\eta|x|}|). \ea

Now define \ba\label{xi} \xi(t)&:=\sup_{x, \, 0\le s\le t}
\Big((|\Delta u^n|+|\Delta u^n_x|)(\theta+\psi_1+\psi_2)^{-1}(x,s)
+|\Delta\dot\delta^n(s)(1+s)| \Big). \ea

Let $r':=|(\Delta\delta^{n-1}, \Delta\delta_*^{n-1})|_*=
|\Delta\delta^{n-1}|_{B^1} + M|\Delta\delta_*^{n-1}|,$  $r'$
sufficiently small. From (\ref{varEbounds}) and smallness of
$|u^n|_{H^5}$ and $|\hat u^n|_{H^5}$, and the fact that $r,r'<<
1$, we obtain
\ba |\Delta u^n(\cdot, t)|_{H^4}\le
C(r'+\xi(t))(1+t)^{-\frac14},\ea which gives us  a bound on
$L^\infty$-norm of $\Delta u^n_{xxx}$, providing us, therefore,
with the bounds,
 \ba\label{qrsvar}
|\Delta Q +\Delta R|&\le C(r\xi(t)+rr')(\Psi+ \Phi_1), \\
|\Delta Q_y +\Delta R_y|, \, |\Delta Q _{yy} +\Delta R_{yy}|&\le
C(r\xi(t)+rr') \Upsilon \\
|\Delta S|,\, |\Delta S_y|&\le C(r\xi(t)+rr') \Phi_2, \ea
Also, \eqref{sourcebds} and \eqref{derivsourcebds}  hold with
$c(\zeta^2 + 4C^2E_0^2)$ replaced with $Cr^2$. We use
(\ref{Hderivbds}), (\ref{ederivbds}) and (\ref{Gderivbds}) to
obtain
\ba\label{deltae}\Delta e^{n-1} \sim e\Delta\des \le r' e, \ea
 and
similar appropriate bounds for $\Delta H^{n-1}, \Delta \tilde
G^{n-1}$ and their derivatives (of course, $e$ in (\ref{deltae})
is defined at a point between $\delta_*^{n-1}$ and
$\hat\delta_*^{n-1}$). Next, using lemmas
\ref{iniconvolutions}--\ref{Hnonlinear} in a procedure parallel to
the one used in the proof of Lemma \ref{uexistence}, we obtain
$$\xi (t) \le C(r' + r\xi(t))$$
from which we conclude that $$\xi(t)\le \frac{Cr'}{1-Cr}$$ with
constant $C$ independent of $r$ and $r'$.  Now, replacing $\xi$ in
(\ref{qrsvar}) with this bound,  we plug back
 the result into (\ref{vardotdeltan}) and into
the similar formula for $\Delta \delta^n$. Notice that, with the
exception of the first two terms, the other term in
\eqref{vardotdeltan} have quadratic terms in their source term, so
giving us small enough bounds. Hence, using one again lemmas
\ref{iniconvolutions}--\ref{Hnonlinear},  we obtain
\ba \label{vardotdeltabd} |\Delta\dot\delta^n|\le (CE_0r'
+\frac{C}{M}r'+ Crr')(1+t)^{-1},\ea of which the two first terms
in the right hand side come from the first two terms of
\eqref{vardotdeltan}. Similarly we obtain
\ba\label{vardeltabd} |\Delta\delta^n|\le (CE_0 +\frac{C}{M}+
Cr)r'(1+t)^{-\frac12}\ea
We notice that $(CE_0 +\frac{C}{M}+ Cr)$ can be made arbitrarily
small, provided that $E_0$, $r$ are small enough and $M$ is large
enough. Next, we use \eqref{vardeltasn} to bound $\Delta\des^n$,
using basically the same method used in
(\ref{deslin})--(\ref{deslina}), and therefore obtaining
 \ba M|\Delta\delta_*^n| \le (CE_0+ Cr)r'.\ea
This, together with \eqref{vardotdeltabd} and \eqref{vardeltabd},
gives us (\ref{lipbds}) with $\alpha < 1$, finishing the proof of
the (main) Theorem \ref{nonlin}.

\end{proof}
%CHANGED: grammar only-K
\begin{rem}\textup{ In order to control the 
%\begin{rem}\textup{ In order to control 
$H^4$ norm of the variational problem, as in
(\ref{varEbounds}), we indeed need regularity $C^5$
%(\ref{varEbounds}), we indeed need the regularity of class $C^5$
for the coefficients in hypothesis (H0), 
%for one derivative is being lost in the variational energy estimates.}
since one derivative is lost in variational energy estimate \eqref{Ebounds}.}
%ENDCHANGED: 
\end{rem}

\br\label{strat}
\textup{
HERE
}
\er


\begin{thebibliography}{GMWZ2}
{\footnotesize

\bibitem[AMPZ1]{AMPZ1} A. Azevedo, D. Marchesin, B. Plohr and K. Zumbrun,
{\it Nonuniqueness of Riemann solutions},
Z. Angew. Math. Phys. 47 (1996), no. 6, 977--998.

\bibitem[Br1]{Br1} L. Q. Brin, {\it Numerical testing of the stability of viscous
shock waves}, Ph.D. dissertation, Indiana University, May 1998.

\bibitem[Br2]{Br2} L. Q. Brin,
{\it Numerical testing of the stability of viscous shock waves,}
Math. Comp. 70 (2001) 235, 1071--1088.

\bibitem[Br3]{Br3} L. Brin,
{\it Numerical testing of the stability of viscous shock waves},
Doctoral thesis, Indiana University (1998).

\bibitem[BrZ]{BrZ} L. Brin and K. Zumbrun,
{\it Analytically varying eigenvectors and the stability of viscous
shock waves}. Seventh Workshop on Partial Differential Equations, Part I (Rio de Janeiro, 2001).
Mat. Contemp. 22 (2002), 19--32.

\bibitem[BDG]{BDG}
T. Bridges, G. Derks, and G. Gottwald,
{\it Stability and instability of solitary waves of the fifth-order
KdV equation: a numerical framework,}
Phys. D  172  (2002),  no. 1-4, 190--216.

\bibitem[Fre]{Fre} H. Freist\"uhler,
{\it Some results on the stability
of non-classical shock waves,} J. Partial Diff. Eqs. 11 (1998), 23-38.

\bibitem[FreS]{FreS} H. Freist\"uhler and P. Szmolyan,
{\it Spectral stability of small shock waves,}
Arch. Ration. Mech. Anal. 164 (2002) 287--309.

\bibitem[FreZ]{FreZ} H. Freist\"uhler and K. Zumbrun,
{\it Examples of unstable viscous shock waves,}
unpublished note, Institut f\"ur Mathematik, RWTH Aachen, February 1998.

\bibitem[GZ]{GZ} R. Gardner and K. Zumbrun,
{\it The Gap Lemma and geometric criteria for instability
of viscous shock profiles},
Comm. Pure Appl.  Math. 51 (1998), no. 7, 797--855.

\bibitem[HZ]{HZ} P. Howard and K. Zumbrun,
\emph{Stability of undercompressive viscous shock waves}, in
press, J. Differential Equations  225  (2006),  no. 1, 308--360.

\bibitem[HR]{HR} P. Howard and M. Raoofi,
{\it Pointwise asymptotic behavior of perturbed viscous shock
profiles},
Adv. Differential Equations  11  (2006),  no. 9, 1031--1080.

\bibitem[HRZ]{HRZ} P. Howard, M. Raoofi, and K. Zumbrun,
\emph{Sharp pointwise bounds for perturbed shock waves}, J.
Hyperbolic Differ. Equ. 3 (2006), no. 2, 297--374.

\bibitem[HuZ]{HuZ} J. Humpherys and K. Zumbrun,
{\it Spectral stability of small amplitude
shock profiles for dissipative symmetric hyperbolic--parabolic systems,}
Z. Angew. Math. Phys. 53 (2002) 20--34.

\bibitem[Kaw]{Kaw} S. Kawashima,
\emph{Systems of a hyperbolic--parabolic composite type,
 with applications to the equations of magnetohydrodynamics},
thesis, Kyoto University (1983).

%\bibitem{KM} S. Kawashima and A. Matsumura,
%{\it Asymptotic stability of traveling wave solutions of systems
%for one-dimensional gas motion,}
%Comm. Math. Phys. 101 (1985), no. 1, 97--127.
%
%\bibitem{KMN} S. Kawashima, A. Matsumura, and K. Nishihara,
%{\it Asymptotic behavior of solutions for the equations of a
%viscous heat-conductive gas,}
%Proc.  Japan Acad. Ser. A Math. Sci. 62 (1986), no. 7, 249--252.

%\bibitem{La} P.D. Lax,
%{\it Hyperbolic systems of conservation laws and the mathematical
%theory of shock waves},
%Conference Board of the Mathematical Sciences Regional Conference
%Series in Applied Mathematics, No. 11.
%Society for Industrial and Applied Mathematics,
%Philadelphia, Pa., 1973. v+48 pp.


%\bibitem{LZ.1} T.P. Liu and K. Zumbrun,
%{\it Nonlinear stability of an undercompressive shock for complex
%Burgers equation,} Comm. Math. Phys. 168 (1995), no. 1, 163--186.
%
%\bibitem[LZe]{LZe} T.-P. Liu and Y. Zeng,
%{\it Large time behavior of solutions for general
%quasilinear hyperbolic--parabolic systems of conservation laws}.
%AMS memoirs 599 (1997).

\bibitem[LZ2]{LZ2} T.P. Liu and K. Zumbrun,
{\it On nonlinear stability of general undercompressive viscous shock waves,}
Comm.  Math. Phys. 174 (1995), no. 2, 319--345.

\bibitem[LRTZ]{LRTZ} G. Lyng, M. Raoofi, B. Texier, and K. Zumbrun,
\emph{Pointwise Green Function Bounds and stability of combustion waves},
to appear, J. Diff. Eq. (2007).

\bibitem[M1]{M1} A. Majda,
{\it The stability of multi-dimensional shock fronts -- a
new problem for linear hyperbolic equations,}
Mem. Amer. Math. Soc. 275 (1983).

\bibitem[M2]{M2} A. Majda,
{\it The existence of multi-dimensional shock fronts,}
Mem. Amer. Math. Soc. 281 (1983).

\bibitem[M3]{M3} A. Majda,
{\it Compressible fluid flow and systems of conservation laws in several
space variables,} Springer-Verlag, New York (1984), viii+ 159 pp.

%
\bibitem[MP]{MP} A. Majda and R. Pego, \textit{
Stable viscosity matrices for systems of conservation laws}, J.
Diff. Eqs. 56 (1985) 229--262.
%

%\bibitem{MN} A. Matsumura and K. Nishihara,
%{\it On the stability of travelling wave solutions of a one-dimensional
%model system for compressible viscous gas,}
%Japan J. Appl. Math. 2 (1985), no. 1, 17--25.
%
\bibitem[MZ1]{MZ1} C. Mascia and K. Zumbrun,
{\it Pointwise Green's function bounds and stability of relaxation shocks},
Indiana Univ. Math. J. 51 (2002), no. 4, 773--904.

\bibitem[MZ2]{MZ2} C. Mascia and K. Zumbrun,
{\it Stability of small-amplitude shock profiles of symmetric
hyperbolic-parabolic systems,}
Comm. Pure Appl. Math.  57  (2004),  no. 7, 841--876.

\bibitem[MZ3]{MZ3} C. Mascia and K. Zumbrun,
{\it Pointwise Green's function bounds for shock profiles
with degenerate viscosity,}
Arch. Rational Mech. Anal. 169  (2003),  no. 3, 177--263.

\bibitem[MZ4]{MZ4} C. Mascia and K. Zumbrun,
{\it Stability of large-amplitude shock profiles of hyperbolic--parabolic
systems},
Arch. Rational Mech. Anal. 172  (2004),  no. 1, 93--131.

\bibitem[MZ5]{MZ5} C. Mascia and K. Zumbrun,
{\it Stability of large-amplitude shock profiles of general relaxation
systems},
SIAM J. Math. Anal.  37  (2005),  no. 3, 889--913.

\bibitem[PZ]{PZ} R. Plaza and K. Zumbrun,
{\it An Evans function approach to spectral stability
of small-amplitude viscous shock profiles,}
J. Disc. and Cont. Dyn. Sys. 10. (2004), 885-924.

\bibitem[Ra]{Ra} M.-R. Raoofi, {\it $L^p$-asymptotic behavior
of perturbed viscous shock profiles,}
 J. Hyperbolic Differ. Equ.  2  (2005),  no. 3, 595--644.

%SOME OLD UC REFS...
%\bibitem{Sl.1} M. Slemrod,
%{\it The viscosity-capillarity approach to phase transitions,} in:
%{\it PDEs and continuum models of phase transitions} (Nice, 1988), 201--206,
%Lecture Notes in Phys., 344,
%Springer, Berlin-New York, 1989.
%
%\bibitem{Sl.2} M. Slemrod,
%{\it The vanishing viscosity-capillarity approach
%to the Riemann problem for a van der Waals fluid,}
%in: {\it Nonclassical continuum mechanics} (Durham, 1986), 325--335,
%London Math. Soc. Lecture Note Ser., 122,
%Cambridge Univ. Press, Cambridge-New York, 1987.
%
%\bibitem{Sl.3} M. Slemrod,
%{\it A limiting ``viscosity'' approach to the
%Riemann problem for materials exhibiting change of phase,}
%Arch. Rational Mech. Anal. 105 (1989), no. 4, 327--365.
%
%\bibitem{Sl.4} M. Slemrod,
%{\it Dynamic phase transitions in a van der Waals fluid,}
%J. Differential Equations 52 (1984), no. 1, 1--23.
%
%\bibitem{Sl.5} M. Slemrod,
%{\it Admissibility criteria for propagating phase
%boundaries in a van der Waals fluid,}
%Arch. Rational Mech. Anal. 81 (1983), no. 4, 301--315.
 %
%\bibitem{Sm} J. Smoller,
%{\it Shock waves and reaction--diffusion equations,}
%Second edition, Grundlehren der Mathematischen Wissenschaften,
%{Fundamental Principles of Mathematical Sciences}, 258.
%Springer-Verlag, New York, 1994. xxiv+632 pp. ISBN: 0-387-94259-9.
%

\bibitem[TZ1]{TZ1} B. Texier and K. Zumbrun,
{\it Galloping instability of viscous shock waves,}
Preprint (2006); available http://arxiv.org/abs/math.AP/0609331.
%TODO: update?

\bibitem[TZ2]{TZ2} B. Texier and K. Zumbrun,
{\it Hopf bifurcation of viscous shock waves in
compressible gas- and magnetohydrodynamics,}
preprint (2006); available http://arxiv.org/abs/math.AP/0612044.
%TODO: update?


\bibitem[ZH]{ZH} K. Zumbrun and P. Howard,
{\it Pointwise semigroup methods and stability of viscous shock waves,}
Indiana Mathematics Journal V47 (1998) no. 4, 741--871.

%\bibitem{Z.1} K. Zumbrun,  {\it Stability of viscous shock waves},
%Lecture Notes, Indiana University (1998).

%\bibitem{Z2}[Z1] K. Zumbrun, {\it Refined Wave--tracking and Nonlinear
%Stability of Viscous Lax Shocks}, Methods Appl. Anal.  7 (2000) 747--768.

\bibitem[Z1]{Z3} K. Zumbrun, {\it Multidimensional stability of
planar viscous shock waves},
Advances in the theory of shock waves, 307--516,
Progr. Nonlinear Differential Equations Appl., 47, Birkhäuser Boston,
Boston, MA, 2001.

\bibitem[Z2]{Z4} K. Zumbrun, {\it Stability of large-amplitude shock
waves for compressible Navier--Stokes equations,}
to appear, Handbook of Fluid Mechanics, Elsevier (2004).

\bibitem[Z3]{Z5} K. Zumbrun,
{\it Planar stability criteria for
viscous shock waves of systems with real viscosity,}
to appear, CIME lecture notes series; preprint (2004).

%\bibitem{Z.6} K. Zumbrun, {\it Dynamical stability of phase transitions
%in the $p$-system with viscosity-capillarity,}
%SIAM J. Appl. Math. 60 (2000) 1913--1924.

\bibitem[ZS]{ZS} K. Zumbrun and D. Serre,
{\it Viscous and inviscid stability of multidimensional
planar shock fronts,} Indiana Univ. Math. J. 48 (1999) 937--992.

}
\end{thebibliography}
\end{document}